%
%
%

\input amstex  

\input amssym
\input amssym.def

\magnification 1200
\loadmsbm
\parindent 0 cm

\define\nl{\bigskip\item{}}
\define\snl{\smallskip\item{}}
\define\inspr #1{\parindent=20pt\bigskip\bf\item{#1}}
\define\iinspr #1{\parindent=27pt\bigskip\bf\item{#1}}
\define\einspr{\parindent=0cm\bigskip}

\define\tussenen{\quad\qquad\text{and}\qquad\quad}

\define\ot{\otimes}

\define\tl{\triangleleft}

\centerline{\bf Separability Idempotents and Multiplier Algebras}
\bigskip
\centerline{\it  Alfons Van Daele \rm $^{(*)}$}
\bigskip\bigskip
{\bf Abstract} 
\nl 
Consider two non-degenerate algebras $B$ and $C$ over the complex numbers. We study a certain class of idempotent elements $E$ in the multiplier algebra $M(B\ot C)$, called further {\it separability idempotents}. The main requirement is that there are non-degenerate anti-homomorphisms $S_B$ from $B$ to the multiplier algebra $M(C)$ and $S_C$ from $C$ to the multiplier algebra $M(B)$ of $B$ satisfying 
$$E(b\ot 1)=E(1\ot S_B(b))
\qquad\qquad\text{and}\qquad\qquad
(1\ot c)E=(S_C(c)\ot 1)E.$$
for all $b\in B$ and $c\in C$. There exist linear functionals $\varphi_B$ and $\varphi_C$ on $B$ and $C$ respectively, satisfying and characterized by  
$$(\varphi_B\ot \iota)E=1
\qquad\qquad\text{and}\qquad\qquad
(\iota\ot\varphi_C)E=1.$$
Here $\iota$ denotes the identity map, both on $B$ and on $C$. These formulas hold in the multiplier algebras $M(C)$ and $M(B)$ respectively. 
\snl
We also treat the notion of regularity for such a separability idempotent. We call $E$ {\it regular} if it is a separability idempotent also when considered in $M(B^{\text{op}}\ot C^{\text{op}})$. Here we use $B^{\text{op}}$ and $C^{\text{op}}$ for the algebras that are obtained from $B$ and $C$ respectively by reversing the product. We show that $E$ is regular if and only if the anti-homomorphisms $S_B$ and $S_C$ have range in $C$ and $B$ respectively and are actually anti-isomorphisms.
Regularity of a separability idempotent is automatic if the underlying algebras are finite-dimensional. The same is true if they are abelian or if they are $^*$-algebras and if $E$ is self-adjoint. 
\snl
In the case of a regular separability idempotent, we are able to determine the structure of the underlying algebras $B$ and $C$. We also obtain some partial results for the non-regular case. We discuss further problems at the end of the paper.
\snl
Separability idempotents naturally appear in the theory of weak multiplier Hopf algebras and quantum groups. If $(A,\Delta)$ is a weak multiplier Hopf algebra, then the canonical idempotent $E$ is a separability idempotent in $M(B\ot C)$ where $B$ is the source algebra and $C$ the target algebra. The separability idempotent in this case is regular if the underlying weak multiplier Hopf algebra is regular. 
Another example, of a different nature, is coming with a discrete quantum group $(A,\Delta)$. Here we take $B=C=A$ and $E=\Delta(h)$ where $h$ is the normalized cointegral. The anti-homomorphisms coincide with the original antipode and the linear functionals $\varphi_B$ and $\varphi_C$ are the right and the left integral. 
\nl 
Date: {\it 22 September  2015} 
\bigskip
\hrule
\bigskip
\parindent 0.7 cm
\item{(*)} Department of Mathematics, University of Leuven, Celestijnenlaan 200B,\newline
B-3001 Heverlee, Belgium. E-mail: Alfons.VanDaele\@wis.kuleuven.be

\parindent 0 cm

\newpage

\bf 0. Introduction \rm
\nl
Consider the algebra $M_n(\Bbb C)$ of $n\times n$ complex matrices and denote it by $A$. It is an algebra over $\Bbb C$ with an identity. Let $(e_{ij})$ be a set of matrix units in $A$. Define $E\in A\ot A$ by
$$E=\frac1n \sum_{i,j=1}^n e_{ij}\ot e_{ij}.$$
Then $E$ has the following properties.
\snl
First it is {\it an idempotent}, that is $E^2=E$ in the algebra $A\ot A$. If we consider $A$ with its natural $^*$-algebra structure, then $E$ is {\it self-adjoint}. Next, the legs of $E$ are all of $A$ in the following sense. The left leg of $E$ is defined as the smallest subspace $V$ of $A$ so that $E\in V\ot A$. Then we must have $V=A$. Similarly for the right leg. 
\snl
Consider the linear map $S:A\to A$ defined by $S(e_{ij})=e_{ji}$. Then $S$ is an anti-isomorphism of $A$ and it satisfies
$$E(a\ot 1)=E(1\ot S(a)) \qquad\qquad\text{and}\qquad\qquad (1\ot a)E=(S(a)\ot 1)E$$
for all $a\in A$. Here $S$ is involutive in the sense that $S^2=\iota$ (where we use $\iota$ for the identity map on $A$). It also satisfies $S(a^*)=S(a)^*$ for all $a$ when we consider $A$ with its natural $^*$-algebra structure. We have that 
$$m(S\ot \iota)E=1 \qquad\qquad\text{and}\qquad\qquad m(\iota\ot S)E=1$$
where again $\iota$ is the identity map and where $m$ stands for the multiplication map on $A$, seen as a linear map from $A\ot A$ to $A$. Finally, if $\text{Tr}$ is the trace on $A$, normalized so that $\text{Tr}(1)=n$, and if $\varphi$ is defined as $n\text{Tr}$, we have
$$(\iota\ot\varphi)E=1 \qquad\qquad\text{and}\qquad\qquad (\varphi\ot \iota)E=1.$$
This means that $A$ is a {\it separable algebra} and that $E$ is a {\it separability idempotent} in $A\ot A$ with the extra property that the legs of $E$ are all of $A$. It is therefore {\it separable Frobenius} in a sense that will be explained.
\snl
In this paper, we will consider the case of algebras that are {\it not necessarily finite-dimensio\-nal} and are also {\it not required to be unital}. We will consider separability idempotents in a sense that will be close to the above description. We will show how the known concept in the case of an algebra with identity relates with our approach to the more general case.
\snl
We have started to study this concept as a consequence of our work on weak multiplier Hopf algebras where it naturally pops up. This has influenced our approach. We will explain this and the relation with weak multiplier Hopf algebras further.
\nl
\it Separability idempotents \rm
\nl
Our {\it starting point} is a pair of non-degenerate algebras $B$ and $C$ over the complex numbers and an idempotent element $E$ in the multiplier algebra $M(B\ot C)$ of the tensor product algebra $B\ot C$, satisfying certain requirements as above. Some of the properties described in the simple case of $M_n(\Bbb C)$ will become axioms while some others will be proven from the axioms. As it turns out however, not all of the properties will remain true and some need to be weakened. 
\snl
The main differences are due to the fact that we do not assume the algebras to be finite-dimensional and unital. We will have the equivalent of the anti-isomorphism $S$ above. As we have two algebras $B$ and $C$, and not a single one as above, we will have anti-homomorphisms $S_B$ and $S_C$ on $B$ and $C$ respectively. In fact, the separability idempotents we study in this paper should perhaps be called {\it Frobenius separability idempotents} because we assume that they are full in the sense that their left and right legs are all of $B$ and $C$ respectively, just like in the motivating example with $M_n(\Bbb C)$ above.  
\snl
We refer to Section 1 of this paper for a precise definition of a {\it separability idempotent} and the general results, in particular about these {\it antipodal maps} on $B$ and $C$.    
\nl
When writing this paper, we have been {\it motivated by two special cases}.
\snl
The first most important case comes from a {\it weak multiplier Hopf algebra} $(A,\Delta)$. Here the algebras $B$ and $C$ are the images $\varepsilon_s(A)$ and $\varepsilon_t(A)$ of the source and target maps $\varepsilon_s$ and $\varepsilon_t$ and $E$ is the canonical idempotent playing the role of $\Delta(1)$ for this weak multiplier Hopf algebra. The antipode $S$ gives the anti-homomorphisms $S_B$ and $S_C$. If the weak multiplier Hopf algebra is regular, then so it the separability idempotent. This is characterized by the fact that $S$ is a bijective anti-isomorphism of $A$.
 It is this case that motivated us to start with the investigation of separability idempotents as studied in this paper. We have e.g.\ used the theory in [VD-W4.v2] to construct examples of weak multiplier Hopf algebras. The paper [VD-W4.v2] is an improved version of the older one [VD-W4.v1]. We refer to Section 3 for this example.
\snl
For another case, we consider a {\it discrete quantum group}. By this we mean a multiplier Hopf algebra $(A,\Delta)$ with a (left) cointegral $h$ so that $\varepsilon(h)=1$ where $\varepsilon$ is the counit of $(A,\Delta)$. Then we consider $\Delta(h)$ as sitting inside the multiplier algebra $M(A\ot A)$. This will be a (regular) separability idempotent in the sense of this paper and the antipode $S$ will play the role of the anti-isomorphisms. It is known that $S$ not necessarily satisfies $S^2=\iota$ and that there are multiplier Hopf $^*$-algebras that are discrete quantum groups where $S$ is not a $^*$-map. The trace $\tau$ in the case of $M_n(\Bbb C)$ will be replaced by the left and the right integral on the discrete quantum group. We again refer to Section 3 
where this is explained in detail (and where also references are given). 
\nl
\it Content of the paper \rm 
\nl
In {\it Section} 1 we first give the main definitions and the main results. As mentioned already, we start with a pair of non-degenerate algebras $B$ and $C$ over the complex numbers and an idempotent $E$ in the multiplier algebra $M(B\ot C)$. The first main requirement is that the left and the right legs of $E$ are all of $B$ and $C$ respectively. To be able to formulate this condition correctly, we need a technical requirement first. The next main condition is the existence of non-degenerate anti-homomorphisms $S_B:B\to M(C)$ and $S_C:C\to M(B)$ satisfying 
$$E(b\ot 1)=E(1\ot S_B(b))
\quad\qquad\text{and}\qquad\quad
(1\ot c)E=(S_C(c)\ot 1)E$$
for $b\in B$ and $c\in C$. Here $M(B)$ and $M(C)$ are the multiplier algebras of $B$ and $C$ respectively. It is a consequence of the first requirement about the legs of $E$ that these anti-homomorphisms are unique if they exist and that they are injective. These anti-homomorphisms play a crucial role in the further treatment. They are referred to as the antipodal maps associated with the separability idempotent.
\snl
One of the results we obtain is the existence of local units for the underlying algebras $B$ and $C$.
\snl
We also find that there are unique linear functionals $\varphi_B$ on $B$ and $\varphi_C$ on $C$ satisfying and characterized by
$$(\varphi_B\ot\iota)E=1
\tussenen
(\iota\ot\varphi_C)E=1.$$
The formulas make sense in the multiplier algebras. We begin the study of these {\it distinguished} linear functionals in this section but we obtain nicer results in Section 2.
\snl
In the semi-regular case, the anti-homomorphisms $S_B$ and $S_C$ have range in $C$ and $B$ respectively while in the regular case, they are moreover anti-isomorphisms.
These cases are studied in {\it Section} 2. We give different characterizations. 
\snl
More and nicer results can be proven for regular separability idempotents. We show e.g.\ that the distinguished linear functionals $\varphi_B$ and $\varphi_C$ are faithful and admit modular automorphisms.
\snl
In many cases, regularity is automatic as we show in {\it Section} 3. This is the case when the underlying algebras are finite-dimensional. Also when the algebras $B$ and $C$ are abelian, regularity is automatic. Finally, the involutive case is considered as well. Here we assume that the algebras $B$ and $C$ are $^*$-algebras and that $E$ is a self-adjoint idempotent in $M(B\ot C)$. Again, in this situation, $E$ is always regular.
\snl
In {\it Section} 4 we study modules over the base algebras $B$ and $C$ in the case when there is a separability idempotent $E\in M(B\ot C)$. If the separability idempotent is regular, we prove that unital modules split. In general and in the semi-regular case, there are  some similar results, but they are more subtle.
\snl
In the regular case however, where we do have nice splitting results, we can apply this and obtain that the base algebras have to be semi-simple, just as in the finite-dimensional case and we obtain a complete characterization of the regular separability idempotents.   
\snl
Finally, in {\it Section} 5, we draw some conclusions and we discuss possible further research. The regular case is completely understood as we have seen in Section 4. However, if $E$ is not regular, many questions are still open. This is also the case for the semi-regular case. In particular, unfortunately examples of non-regular separability idempotents are not (yet) available.

\newpage
\it Conventions and notations \rm
\nl
We only work with algebras $A$ over $\Bbb C$ (although we believe that this is not essential and algebras over other fields can also be considered). We do not assume that they are unital but we need that the product is non-degenerate. By this we mean that the product as a bilinear map is non-degenerate. As it will turn out, the algebras will have local units. In particular, they are idempotent. This means that any element is a sum of elements $ab$ where both $a$ and $b$ are elements in the algebra. We write this as $A=A^2$.
\snl
Recall that an algebra $A$ is said to have {\it left local units} if for any finite set $\{a_1,a_2,\dots,a_n\}$, there exists an element $e\in A$ so that $ea_i=a_i$ for all $i$. In [Ve], it is shown to be sufficient only to require that for any element $a$ in $A$ there is an element $e\in A$ satisfying $ea=a$. This result is obtained by induction. Similarly, right local units are defined, with a similar property. The algebra $A$ has two sided local units, or simply local units, if for any finite set $\{a_1,a_2,\dots,a_n\}$, there exists an element $e\in A$ so that $ea_i=a_i$ and $a_ie=a_i$ for all $i$. In [Ve] it is also proven that local units exists if left local units and right local units exist. Indeed, if we have $e',e''$ in $A$ satisfying $e'a_i=a_i$ and $a_ie''=a_i$ for all $i$, then $ea_i=a_i$ and $a_ie=a_i$ when $e$ is defined as $e'+e''-e''e'$
\snl
When $A$ is a non-degenerate algebra, we use $M(A)$ for the multiplier algebra of $A$. When $m$ is in $M(A)$, then by definition we can define $am$ and $mb$ in $A$ for all $a,b\in A$ and we have $(am)b=a(mb)$. The algebra $A$ sits in $M(A)$ as an essential two-sided ideal and $M(A)$ is the largest algebra with identity having this property. If already $A$ has an identity, then $M(A)=A$. If not then $M(A)$ will be strictly larger than $A$.
\snl
If $A$ and $B$ are non-degenerate algebras and if $\alpha:A\to M(B)$ is a homomorphism, it is called {\it non-degenerate} if $\alpha(A)B=B\alpha(A)=B$. In that case, it has a unique extension to a unital homomorphism from $M(A)$ to $M(B)$. The extension is still denoted by $\alpha$. A similar result is true for anti-homomorphisms.
\snl
For a non-degenerate algebra $A$, we consider $A\ot A$, the tensor product of $A$ with itself. It is again a non-degenerate algebra and we can consider the multiplier algebra $M(A\ot A)$. The same is true for a multiple tensor product. We use $\zeta$ for the flip map on $A\ot A$, as well as for its natural extension to $M(A\ot A)$.
\snl
When $A$ is an algebra, we denote by $A^{\text{op}}$ the algebra obtained from $A$ by reversing the product. 
\snl
We use $1$ for the identity in any of these  multiplier algebras. On the other hand, we 
use $\iota$ for the identity map on $A$ (or other spaces).
\snl
A linear functional $f$ on $A$ is called {\it faithful} if the bilinear map $A\times A\to \Bbb C$, mapping $(a,b)$ to $f(ab)$ is non-degenerate. So, given $a\in A$, we have that $a=0$ if either $f(ab)=0$ for all $b$ or $f(ba)=0$ for all $b$. If only one of the conditions is satisfied, we call $f$ left, respectively right faithful. More precisely, $f$ is left faithful if $f(ba)=0$ for all $b$ implies $a=0$ and it is called right faithful if $f(ba)=0$ for all $a$ implies $b=0$.
In the finite-dimensional case, one implies the other and so the notion is not different from faithfulness. In fact, in this case, the existence of a faithful functional is equivalent with the algebra being Frobenius. 
\snl
When we have a $^*$-algebra and when $f$ is positive (i.e\ $f(a^*a)\geq 0$ for all $a$), then $f$ is faithful if and only if $f(a^*a)=0$ implies $a=0$.
\snl
Some of our examples will involve matrix algebras. We will use $M_n(\Bbb C)$ or simply $M_n$ for the ($^*$)-algebra of all complex $n\times n$ matrices where $n=1,2,3,\cdots$
\snl
Finally, when $V$ is any vector space, we will use $V'$ for its linear dual space and evaluation of a functional $\omega\in V'$ in an element $v\in V$ will often be written as $\langle v,\omega\rangle$ or $\langle \omega,v\rangle$.
So we will in particular not bother about the order in this case.
\nl
\it Basic references \rm
\nl
For the theory of separability for algebras with identity, we refer to [Pe] and [W]. In this paper, we also refer often to our work on weak multiplier Hopf algebras. See [VD-W2] for an introduction to the subject (with a motivation of the axioms), [VD-W3] for the main theory and [VD-W4.v2] for the results on the source and target maps and the source and target algebras. Especially this last paper is intimately related with this work on separability idempotents. Finally we are also using discrete quantum groups to illustrate some of our results. Here we refer to [VD-Z].
\snl
Occasionally we will refer to [VD3.v1]. It is a first version of this paper, still available on the arXiv. In that paper, we only considered what we call here {\it regular} separability idempotents. Even in that special case, the present version contains better results. The two versions are substantially different, but as they treat the same objects, we have decided to present this paper as a new version of the first one. Recall also that [VD-W4.v1] is an older version of [VD-W4.v2] that is based on the earlier version of this paper [VD3.v1].
\nl
\bf Acknowledgments \rm
\nl
I would like to thank Joost Vercruysse for discussions about this subject and for introducing me to the existing theory on separability. I am also indebted to Leonid Vainerman for a more recent discussion about the relation of my work with the existing theory of separability for algebras with identity.
\nl\nl

\bf 1. Separability idempotents \rm
\nl
In this section we start with the definition of a separability idempotent and we develop the theory. In the next section, we will consider the important special case of a {\it regular} separability idempotent. We  illustrate a few aspects with simple examples in Section 3 where we first consider the well-known case of finite-dimensional algebras with identity. We then also give more and less trivial examples, as well as some other special cases. 
\nl
\it The main definition \rm
\nl
We start with two {\it non-degenerate algebras $B$ and $C$} and an {\it idempotent element $E$ in the multiplier algebra $M(B\ot C)$}. We have natural embeddings
$$B\ot C\subseteq M(B)\ot M(C)\subseteq M(B\ot C).$$
If the algebras $B$ and $C$ have no identity, these inclusions in general are strict, also the last one. These embeddings are implicitly used further in the paper.
\snl
We will {\it only consider idempotents} with the property that
$$E(1\ot c)\in B\ot C \qquad\quad\text{and}\qquad\quad (b\ot 1)E\in B\ot C \tag"(1.1)"$$
for all $b\in B$ and $c\in C$. Here $1$ is the identity in $M(B)$ or in $M(C)$. We will see later that this is quite natural. 
\snl
In the {\it semi-regular case} (see Definition 2.1 in the next section), as well as in the regular case (see Definition 2.4), we also will require  
$$(1\ot c)E\in B\ot C  \qquad\quad\text{and}\qquad\quad E(b\ot 1)\in B\ot C \tag"(1.2)"$$
for all $b\in B$ and $c\in C$.

\snl
Of course, these assumptions are void in the case of algebras with identity as then we have $E\in B\ot C$. But in general, they are genuine assumptions (also the regularity assumption).
\snl
The reader should compare this with the notion of a ({\it regular}) {\it coproduct} on an algebra as defined e.g.\ in Definition 1.3 in [VD-W3]. 
\snl
To formulate the next assumption, we introduce the following notion (see also Definition 1.5 in [VD-W1]). Observe that it needs condition (1.1) above to make sense.

\inspr{1.1} Definition \rm By the {\it left leg} and the {\it right leg} of $E$ we mean the smallest subspaces $V$ of $B$ and $W$ of $C$ respectively satisfying 
$$E(1\ot c)\in V\ot C \qquad\quad\text{and}\qquad\quad (b\ot 1)E\in B\ot W $$
for all $c\in C$ and $b\in B$. We call $E$ {\it full} if the left and right legs of $E$ are all of $B$ and $C$ respectively.
\hfill$\square$\einspr

Again, the reader should compare this with the notion of a {\it full coproduct} as given in Definition 1.4 of [VD-W3]. In the regular case, the left and right legs can also be defined with the 
factors $1\ot c$ and $b\ot 1$ on the other side of $E$. This will give the same subspaces. See similar results in [VD-W1], [VD-W2] and [VD-W3] where related topics are treated.
\nl
We have the following easy consequences.
    
\inspr{1.2} Proposition \rm Let $E$ be full. Then any element in $B$ is a linear combination of elements of the form $(\iota\ot\omega)(E(1\ot c))$ with $c\in C$ and $\omega$ a linear functional on $C$. Similarly, any element in $C$ is a linear combination of elements of the form $(\omega\ot\iota)((b\ot 1)E)$ where $b\in B$ and with $\omega$ a linear functional on $B$. 
\hfill$\square$\einspr
In the regular case, we can also take expressions with $b$ and $c$ on the other side. The proof is easy, see e.g.\ Section 1 in [VD-W1] and in [VD-W2] where  similar results are proven.

\inspr{1.3} Proposition \rm Let $E$ be full. Then we have, for $b\in B$ and $c\in C$, that
$$\align (b\ot 1)E&=0 \qquad \text{implies}\qquad b=0 \\
         E(1\ot c)&=0	\qquad \text{implies}\qquad c=0.
\endalign$$
The same results are true with the factors on the other side. 			 
\hfill$\square$\einspr

This property is an easy consequence of the previous one, using non-degeneracy of the products in $B$ and $C$. Regularity is not needed for the result with the factors on the other side.
\nl
Next, we formulate the most fundamental assumption and we come to the main definition of this section. 

\inspr{1.4} Definition \rm Let $E$ be a full idempotent in $M(B\ot C)$ as defined in Definition 1.1. Assume moreover that there are non-degenerate anti-homomorphisms $S_B:B\to M(C)$ and $S_C:C\to M(B)$ such that 
$$E(b\ot 1)=E(1\ot S_B(b)) 
\qquad\quad\text{and}\quad\qquad
(1\ot c)E=(S_C(c)\ot 1)E. \tag"(1.3)"$$
for all $b\in B$ and $c\in C$. Then we call $E$ a {\it separability idempotent in $M(B\ot C)$}. 	
\hfill$\square$\einspr

Sometimes, we will also refer to $S_B$ and $S_C$ as the {\it antipodal maps}.
\snl
Recall that non-degeneracy of $S_B$ means that elements $S_B(b)c$ with $b\in B$ and $c\in C$ span all of $C$ and the same for elements of the form $cS_B(b)$. Similarly for $S_C$. 
\snl
Also remark that the assumptions in Definition 1.1 will follow automatically from the existence of these non-degenerate anti-homomorphisms. Indeed $E(1\ot S_B(b)c)=E(b\ot c)$ for all $b\in B$ and $c\in C$ and this will belong to $B\ot C$ as $E\in M(B\ot C)$. Because $S_B$ is non-degenerate, this implies that $E(1\ot c)\in B\ot C$ for all $c\in C$. Similarly, the existence of the non-degenerate anti-homomorphism $S_C$ will give that $(b\ot 1)E\in B\ot C$ for all $b\in B$.  
\snl
Remark that the anti-homomorphisms $S_B$ and $S_B$ are unique if they exist. This follows from Proposition 1.3. The same result also gives that the maps $S_B$ and $S_C$ necessarily are injective. In fact, the defining relations (1.3) also imply that the maps are anti-homomorphisms.
\snl
We can deduce from these remarks that the definition is quite natural. 
\snl
Before we continue, we would like to make some more remarks. 

\inspr{1.5} Remark \rm
i) As mentioned already, in the case of unital algebras, some of the above conditions are automatically satisfied. This is the case for the conditions (1.1) and (1.2). However, even for finite-dimensional unital algebras, fullness is {\it not} automatic. 
\snl
ii) Fullness of $E$ is an important property of the separability idempotents as we study them in this paper. In the more familiar finite-dimensional case, it follows that in fact, we consider separability idempotents that are of {\it Frobenius type}. 
\snl
ii) As we will also see later, in Section 2, if $E$ is regular, the conditions will imply that the algebras $B$ and $C$ are anti-isomorphic and that they are separable algebras. However, it is not clear if this is still the case when regularity is not assumed. 
\hfill$\square$\einspr

Later, we will come back to some of these remarks and explain more precisely the use of our terminology and compare it with the existing literature. 
We refer also to Section 5 where we discuss some open problems and possible further research on these objects.
\nl
\it The main properties \rm
\nl
We will now prove several properties of a separability idempotent as defined above. We will see that the proofs are not very difficult. In Section 2, where we treat the regular case, we will find nicer results. These results are illustrated in Section 4. And they are used in [VD-W4.v2] for constructing examples of weak multiplier Hopf algebras. 
\snl
In the next proposition, we consider the multiplication maps $m_C:C\ot C\to C$ and $m_B:B\ot B\to B$ (or rather the obvious extensions to $M(C)\ot C$ and $B\ot M(B)$ respectively).

\inspr{1.6} Proposition \rm 
If $E$ is a separability idempotent in $M(B\ot C)$, then 
$$m_C(S_B\ot\iota)(E(1\ot c))=c
\qquad\qquad\text{and}\qquad\qquad
m_B(\iota\ot S_C)((b\ot 1)E)=b$$ 
for all $c\in C$ and all $b\in B$. 

\snl\bf Proof\rm: 
Take any element $c\in C$ and write $c'$ for $m_C(S_B\ot\iota)(E(1\ot c))$. Then we have
$$E(1\ot c')=E(E(1\ot c))=E(1\ot c)$$
because $E^2=E$. This implies $c'=c$ and it proves the first formula. Similarly for the second formula.
\hfill$\square$\einspr

We think of the formulas above as
$$S_B(E_{(1)})E_{(2)}=1 \qquad\quad\text{and}\qquad\quad E_{(1)}S_C(E_{(2)})=1$$
with the Sweedler type notation $E=E_{(1)}\ot E_{(2)}$. The interpretation is straightforward in the case of unital algebras. In the non-unital case, one still can interpret these formulas in the multiplier algebras $M(C)$ and $M(B)$ (as we have done in the formulation of the results). See also a remark after Proposition 1.9.
\snl
We use this formula in the following result related with the origin of the term 'separability idempotent'.

\inspr{1.7} Proposition \rm 
Assume as before that $E$ is a separability idempotent in $M(B\ot C)$. Consider $B\ot C$ as a right $(B\ot C)$-module with multiplication as the action. Also consider $C$ as a right module over $B\ot C$ with action $x\tl (b\ot c)=S_B(b)xc$ where $S_B$ is the given antipodal from $B$ to $M(C)$. 
Define a linear map $\theta:B\ot C \to C$ by $\theta(b\ot c)=S(b)c$. Then $\theta$ is a module map and the short exact sequence 
$$0\longrightarrow \text{Ker}(m)\longrightarrow B\ot C \overset m\to\longrightarrow C \longrightarrow 0$$
is split.
\snl\bf Proof\rm:
It is straightforward to show that $\theta$ is a module map for the given right module structures on $B\ot C$ and $C$. If we define $\gamma:C\to B\ot C$ by $\gamma(c)=E(1\ot c)$, it follows from the equation $E(b\ot 1)=E(1\ot S_B(b))$ for $b\in B$ that $\gamma$ is a also a module map. And because $m_C (S_B\ot\iota)(E(1\ot c))=c$, we see that $\theta\gamma(c)=c$ for all $c\in C$. This proves the result.
\hfill$\square$\einspr

Remark that we are only using part of the assumptions for the proof of the above result. We just need the anti-homomor\-phism $S_B:B\to M(C)$ and not the other one $S_C:C\to M(B)$. Using the other antipodal map  $S_C$, we can formulate a similar result. 
\snl
Related is the following example. 

\inspr{1.8} Example \rm
Consider again the case where $B$ and $C$ are the algebra $M_n(\Bbb C)$. Now define $E$ in $B\ot C$ as
$$E=\sum_i e_{1i}\ot e_{1i}$$
where we use again matrix units $(e_{ij})$ in $M_n(\Bbb C)$. This element is also an idempotent in $M_n(\Bbb C)\ot M_n(\Bbb C)$. The legs of $E$ are strictly smaller than $M_n(\Bbb C)$. It is true that $E(b\ot 1)=0$ implies $b=0$ and that $E(1\ot c)=0$ implies $c=0$ but this result no longer holds when the elements are on the other side. Similarly we still have $E(b\ot 1)=E(1\ot S_0(b))$ for all $b\in B$ (where $S_0$ is transposition of matrices) but no longer $(1\ot c)E=(S_0(c)\ot 1)E$ for $c\in C$. If we put $F=(S_0\ot \iota)E$ we find that $mF=1$ where $m$ is multiplication as before. Also $F\in C\ot C$ and $(c\ot 1)F=F(1\ot c)$ for all $c$. So, only part of the properties of a separability idempotent remain true. 
\hfill$\square$\einspr

We can consider the following more general version of the example. Take two finite-dimensional unital algebras $B$ and $C$ and assume that $E$ is an idempotent in $B\ot C$. Assume that $S_B$ is an anti-isomorphism from $B$ to $C$ such that $E(b\ot 1)=E(1\ot S_B(b))$. Denote $F=(S_B\ot \iota)E$ and assume that $mF=1$. Then we find $F\in C\ot C$ and $(c\ot 1)F=F(1\ot c)$. Also here we have that for an element $c\in C$ we have $c=0$ if $E(1\ot c)=0$ (because $c=m(F(1\ot c))=m(S_B\ot \iota)(E(1\ot c))$).
In the literature, the element $F$ would be the separability idempotent (in $C\ot C^{\text{op}}$), but it need not be a Frobenius separability idempotent because it is not supposed to be full.
\snl
This illustrates the difference of our notion with the more common one in the literature. In fact, the separability idempotents in this paper are in a sense all of a Frobenius type (because they are assumed to be full).
\nl
We  now look at these elements $(S_B\ot\iota)E$ and $(\iota\ot S_C)E$ for a general separability idempotent as studied in this paper.
\snl
The map $\iota\ot S_C$ is a well-defined and non-degenerate homomorphism from $B\ot C$ to $M(B\ot B^{\text{op}})$. Therefore there is a unique extension, still denoted as $\iota\ot S_C$, to a unital homomorphism from $M(B\ot C)$ to $M(B\ot B^{\text{op}})$. So we can define the element $(\iota\ot S_C)E$ in $M(B\ot B^{\text{op}})$. 
\snl
We collect some properties in the following proposition. And we will give some comments on the result after the proof.

\inspr{1.9} Proposition \rm Denote $(\iota\ot S_C)E$ by $F$. 
We have that $(b\ot 1)F\in B\ot S_C(C)$ for all $b$ in $B$ and $m_BF=1$ in the sense that $m_B((b\ot 1)F)=b$. Again $m_B$ denotes multiplication from $B\ot M(B)$ to $B$. We also have that 
$$(1\ot b)F
\tussenen
F(1\ot b)$$
belong to $B\ot B$ for all $b$. Finally $(x\ot 1)F=F(1\ot x)$ for $x\in S_C(C)$. 
\snl
We have similar results for the element $(S_B\ot\iota)E$.

\snl\bf Proof\rm:
i) The first statement is obvious as $(b\ot 1)E\in B\ot C$ and $S_C:C\to M(B)$. Moreover the formula $m_BF=1$ is a reformulation of the second formula in Proposition 1.6.
\snl
ii) For any $b\in B$ and $c\in C$ we find
$$ 
(1\ot bS_C(c))F
=(1\ot bS_C(c))((\iota\ot S_C)E)
=(1\ot b)(\iota\ot S_C)(E(1\ot c)).
$$
This belongs to $(1\ot B)(B\ot M(B))\subseteq B\ot B$. Because $S_C$ is assumed to be non-degenerate, any element in $B$ is a sum of elements of the form $bS_C(c)$ and so $(1\ot B)F$ is a subset of $B\ot B$.
\snl
iii) To show that also $F(1\ot b)$ belongs to $B\ot B$ is more subtle. We take two copies of $E$ and use the Sweedler type notation $E_{(1)}\ot E_{(2)}$ for one copy and $E^\circ_{(1)}\ot E^\circ_{(2)}$ for the other one. For any $c$ in $C$ we have 
$$E_{(1)}\ot E^\circ_{(1)}\ot E_{(2)} E^\circ_{(2)}c\in B\ot B\ot C.$$ 
This element is equal to
$$E_{(1)}\ot S_C(E_{(2)})E^\circ_{(1)}\ot E^\circ_{(2)}c=(F\ot 1)(1\ot (E(1\ot c)).$$
If we apply a linear functional $\omega$ on the last factor, we find that $F(1\ot b)\in B\ot B$ where $b=(\iota\ot\omega)(E(1\ot c))$. As $E$ is full, it follows that $F(1\ot b)\in B\ot B$ for all $b\in B$.
\snl
iv) Finally we have 
$$(S_C(c)\ot 1)((\iota \ot S_C)E)=(\iota\ot S_C)((1\ot c)E=((\iota\ot S_C)E)(1\ot S_C(c))$$
for all $c\in C$ and this proves the last statement of the proposition.
\hfill $\square$ \einspr 

For some of the arguments to make precise, one needs to cover the formulas in a proper way.
\snl
Because of this result we can also write 
$$E_{(1)}S_C(E_{(2)})b=b
\tussenen
cS_B(E_{(1)})E_{(2)}=c$$
for the result in Proposition 1.6.
\snl
The property that elements like $((\iota\ot S_C)E)(1\ot b)$ belong to $B\ot B$ is a bit surprising. It is somewhat similar to the property we have for (weak) multiplier Hopf algebras $(A,\Delta)$ saying that $((\iota\ot S)\Delta(a))(1\ot b)$ is in $A\ot A$ for all $a,b\in A$ (where $S$ is the antipode). This result holds also in the non-regular case.
\snl
We will reconsider these properties in the next section where we treat regular separability idempotents. We will get better results than the ones above. They will be more similar to the known results in the case of finite-dimensional separable algebras and separability idempotents as studied in this (more restrictive) situation. 
\nl
We now obtain more results on (general) separability idempotents as defined in Definition 1.4 of this section.

\iinspr{1.10} Proposition \rm
If $E$ is a separability idempotent in $M(B\ot C)$, then the algebras $B$ and $C$ have local units. In particular, both algebras are idempotent.

\snl\bf Proof\rm:
We will first show that $b\in bB$ for all $b\in B$. 
\snl
Take any element $b\in B$ and assume that $\omega$ is a linear functional on $B$ so that $\omega$ is $0$ on $bB$. We claim that then $\omega(b)=0$ and this will imply that $b\in bB$. 
\snl
By assumption we have 
$$(\omega\ot\iota)((b\ot 1)E(1\ot c)(y\ot 1))=0$$
for all $c\in C$ and all $y\in M(B)$. This follows from the fact that $E(1\ot c)\in B\ot C$.
On the other hand, we know that $(b\ot 1)E$ belongs to $B\ot C$. Then we can write $(b\ot 1)E$ as a finite sum $\sum p_i\ot q_i$ with $p_i\in B$ and $q_i\in C$ for all $i$. We can also assume that the $(q_i)$ are linearly independent. Then we have $\sum\omega(p_i y)q_ic=0$ for all $c$ and because the product in $C$ is non-degenerate, it follows that $\sum \omega(p_iy)q_i=0$ and so $\omega(p_iy)=0$ for all $i$. This holds for all $y\in M(C)$. Because $S_C(q_i)\in M(C)$ we can conclude that $\omega(p_iS_B(q_i))=0$. If we now take the sum and use that $\sum p_iS_B(q_i)=b$ (as we have shown in the previous proposition), we find $\omega(b)=0$. This proves the claim and it follows that for all $b\in B$ there is an element $e\in B$ so that $be=b$. We know from [Ve] that this implies that $B$ has right local units. 
\snl
Next we show that also $b\in Bb$ for all $b\in B$. This will imply that $B$ also has left local units and again from [Ve], we get that it has two-sided local units (as we mentioned already in the introduction).
\snl
So, let $b\in B$ and now assume that $\omega$ is a linear functional on $B$ so that $\omega$ is $0$ on $Bb$. Take further any two elements $b'$ and $b''$ in $B$. Then by Proposition 1.9,
$$(b''\ot b')(\iota\ot S_C)E\in B\ot B$$ and therefore 
$$(\iota\ot\omega)((b''\ot b')((\iota\ot S_C)E)(1\ot b))=0.\tag"(1.4)"$$
On the other hand, we know that $((\iota\ot S_C)E)(1\ot b)\in B\ot B$ (see Proposition 1.9). We write it as $\sum_i p_i\ot q_i$ and we assume that the $(p_i)$ are linearly independent. We can cancel $b''$ in (1.4) and we get $\omega(b'q_i)=0$ for all $b'\in B$ and all $i$. Now replace $b'$ by $p_i$ and take the sum. Because $\sum_i p_iq_i=m(((\iota\ot S_C)E)(1\ot b))=b$ by Proposition 1.6, it follows that $\omega(b)=0$. Hence $b\in bB$ and this proves the claim.
\snl
In a completely similar way $C$ has local units.
\hfill$\square$\einspr

\nl
\it The distinguished linear functionals $\varphi_B$ and $\varphi_C$\rm
\nl
For every linear functional $\omega$ on $B$ we can define a left multiplier 
$(\omega\ot\iota)E$ of $C$ by $((\omega\ot\iota)E)c=(\omega\ot\iota)(E(1\ot c))$ for all $c\in C$. Similarly we can define a right multiplier $(\iota\ot\omega)E$ of $B$ for any linear functional $\omega$ on $C$ by $b(\iota\ot\omega)E=(\iota\ot\omega)((b\ot 1)E)$.
\snl
With this property in mind, we have the following.

\iinspr{1.11} Proposition \rm
There exist unique linear functionals $\varphi_B$ on $B$ and $\varphi_C$ on $C$ so that
$$(\varphi_B\ot\iota)E=1
\qquad\quad\text{and}\qquad\quad
(\iota\ot\varphi_C)E=1$$

\snl\bf Proof\rm:
i) We first argue that these functionals are unique if they exist. Indeed, if $(\varphi_B\ot \iota)(E(1\ot c))=c$ for all $c\in C$, then $\varphi_B(b')=\omega(c)$ if $b'=(\iota\ot\omega)(E(1\ot c))$ for a given element $c$ and a given linear functional $\omega$ on $C$. Because all elements in $B$ are the linear span of such elements $b'$, uniqueness of $\varphi_B$ will follow.
\snl
ii) The previous argument also suggest how to prove the existence. We have to show that, if 
$$\sum (\iota\ot \omega_i)(E(1\ot c_i))=0 \tag"(1.5)"$$
for a finite number of elements $c_i\in C$ and linear functionals $\omega_i$ on $C$, then we must have $\sum \omega_i(c_i)=0$. If we can show this, then we can define $\varphi_B$ by
$$\varphi_B((\iota\ot \omega)(E(1\ot c))=\omega(c)$$ when $c\in C$ and $\omega$ a linear functional on $C$. Then $\varphi_B$ will satisfy the requirement.
\snl
Now, if we have (1.5), again by the fullness of $E$, we get 
$\sum \omega_i(cc_i))=0$
for all $c$ in $C$. And because $C$ has local units, we can choose $c$ so that $cc_i=c_i$ for all $i$. This implies $\sum \omega_i(c_i)=0$ and so the claim is proven.
\snl
iii) In a completely similar way, we can prove the existence and the uniqueness of the distinguished linear functional $\varphi_C$ on $C$ satisfying $(\iota\ot\varphi_C)E=1$.
\hfill $\square$ \einspr

In the earlier version of this paper (cf.\ [VD3.v1]), we have called these functionals $\varphi_B$ and $\varphi_C$ resp.\ the right and the left integral. This is motivated by the example of a separability idempotent associated with a discrete quantum group (see Proposition 3.12 in Section 3). From another point of view, these functionals are like a counit (see Proposition 2.9 in Section 2). In the sequel we will simply refer to them as the {\it distinguished} linear functionals on $B$ and on $C$.
\snl
Now, we will prove some results about these functionals $\varphi_B$ and $\varphi_C$. First there is the result about faithfulness. 

\iinspr{1.12} Proposition \rm The functional $\varphi_B$ is {\it left} faithful and the functional $\varphi_C$ is {\it right} faithful.

\snl\bf Proof\rm:
Assume that $b\in B$ and that $\varphi_B(b'b)=0$ for all $b'\in B$. This implies that 
$(\varphi_B\ot\iota)(E(b\ot c))=0$ for all $c\in C$. Because
$$(\varphi_B\ot\iota)(E(b\ot c))=(\varphi_B\ot\iota)(E(1\ot S_B(b)c)=S_B(b)c$$
we find that $S_B(b)c=0$ for all $c$. Then $S_B(b)=0$ and because $S_B$ is injective, we get $b=0$.
This proves that $\varphi_B$ is left faithful.
\snl
Similarly, if $c\in C$ and if $\varphi_C(cc')=0$ for all $c'\in C$ we must have $c=0$. So $\varphi_C$ is right faithful.
\hfill $\square$ \einspr

An argument as above does not work to show faithfulness on the other side. Instead, we can show that if $c\in C$ and  $\varphi_B(S_C(c)b')=0$ for all $b'\in B$, then $c=0$. Similarly, if $b\in B$ and if $\varphi_C(c'S_B(b))=0$ for all $c'$, then $b=0$. In the regular case that we consider in the next section, where the maps $S_B$ and $S_C$ have range all of $C$ and $B$ respectively, this will result in faithfulness of both $\varphi_B$ and $\varphi_C$.
\snl
In general, it is not clear if (two-sided) faithfulness of these functionals is always true, even in the semi-regular case.
\nl
Next we consider the so-called {\it weak KMS-property}.

\iinspr{1.13} Proposition \rm The composition $S_C\circ S_B$ is an injective, non-degenerate homomorphism from $B$ to $M(B)$ and 
$$\varphi_B(b'b)=\varphi_B(S_C(S_B(b))b')$$
for all $b,b'\in B$. Similarly, the composition $S_B\circ S_C$ is an injective, non-degenerate homomorphism from $C$ to $M(C)$ and
$$\varphi_C(cc')=\varphi_C(c'S_B(S_C(c)))$$
for all $c,c'\in C$

\snl\bf Proof\rm:
Let $b\in B$. For all $c\in C$ we have
$$(\varphi_B\ot\iota)(E(b\ot c))=(\varphi_B\ot\iota)(E(1\ot S_B(b)c)=S_B(b)c.$$
On the other hand
$$S_B(b)c=(\varphi_B\ot\iota)((1\ot S_B(b))E(1\ot c))=
(\varphi_B\ot\iota)((S_C(S_B(b))\ot 1)E(1\ot c)).$$
It follows that $\varphi(b'b)=\varphi(S_C(S_B(b))b')$ for all $b,b'\in B$. Here $S_C\circ S_B$ is a homomorphism from $B$ to $M(B)$. It is injective because $S_B$ is injective from $B$ to $M(C)$ and because $S_C$ is still injective on $M(C)$. Furthermore we have
$$S_C(S_B(b))S_C(c')b'=S_C(c'S_B(b))b'$$
for all $b,b'\in B$ and $c'\in C$. Now elements of the form $c'S_B(b)$ span all of $C$ and because $S_C$ is non-degenerate it follows that $B$ is spanned by elements of the form $S_C(S_B(b))S_C(c')b'$. This implies in turn that the composition $S_C\circ S_B$ is a non-degenerate homomorphism from $B$ to $M(B)$.
\snl 
Similarly we find that 
$$\varphi_C(cc')=\varphi_C(c'S_B(S_C(c)))$$
for all $c,c'\in C$. And also here the composition $S_B\circ S_C$ is an injective, non-degenerate homomorphism from $C$ to $M(C)$.
\hfill $\square$ \einspr

In the regular case $S_B$ and $S_C$ are anti-isomorphisms and their compositions will be automorphisms (see Proposition 2.8 in the next section). More comments are given there.
\nl
\it Uniqueness properties \rm
\nl 
First we have some kind of uniqueness in the sense that roughly speaking, the separability idempotent is determined by its antipodal maps.

\iinspr{1.14} Proposition \rm
Consider non-degenerate algebras $B$ and $C$. Assume that $E$ and $E'$ are two separability idempotents in $M(B\ot C)$. If the antipodal maps are the same for $E$ and $E'$, then $E=E'$.

\snl\bf Proof\rm:
As before, use $S_B$ and $S_C$ for the antipodal maps. 
\snl
On the one hand we have
$$EE'=E(E'_{(1)}\ot E'_{(2)})=E(1\ot S_B(E'_{(1)})E'_{(2)})=E$$
where we have used the Sweedler notation $E'=E'_{(1)}\ot E'_{(2)}$  and the fact that the anti-homomorphism $S_B$ is the same for $E$ and for $E'$ and Proposition 1.6 for $E'$. On the other hand, we also have
$$EE'=(E_{(1)}\ot E_{(2)})E'=((E_{(1)}S_C(E_{(2)})\ot 1)E'=E'$$
where again we have used the notation $E=E_{(1)}\ot E_{(2)}$ and now the fact that the anti-homomorphism $S_C$ is the same for $E$ and for $E'$ and Proposition 1.6 for $E$. Then $E=E'$.
\hfill$\square$\einspr

To make the argument precise, on can multiply with an element of $C$ on the right in the second factor for proving the first equality and with an element of $B$ on the left in the first factor to prove the second equality.
\snl
This result shows that, if the algebras $B$ and $C$ are given, together with the non-degenerate anti-homomorphisms $S_B:B\to M(C)$ and $S_C:C\to M(B)$, then there is at most one separability idempotent in $M(B\ot C)$ for which these are the associated antipodal maps.
\nl
Next we will consider the idempotent $E'$, defined as $\zeta(S_B\ot S_C)E$ in $M(B\ot C)$. Recall that we use $\zeta$ for the flip map on $B\ot C$ as well as for its extension to the multiplier algebra $M(B\ot C)$. Because the anti-homomorphisms $S_B$ and $S_C$ are assumed to be non-degenerate, the same is true for $S_B\ot S_C$. Then we have a unique extension to an anti-homomorphism on $M(B\ot C)$ to $M(C\ot B)$. Combining it with the flip map $\zeta$ we get $\zeta(S_B\ot S_C)E$ in $M(B\ot C)$. It is obviously again an idempotent.
\snl
A simple argument, as in the proof of the previous proposition, gives that $EE'=E$ and $E'E=E$. Again we need to include the right coverings in order to obtain correct arguments. In this case, we take elements of the form $S_B(b)c$ in the first case and $bS_C(c)$ in the second case
\snl
For this idempotent $E'$ we find the following.

\iinspr{1.15} Proposition \rm We have that $E'$ is a separability idempotent in $M(S_C(C)\ot S_B(B))$. The associated antipodal maps $S'_C$ and $S'_B$ are given by the restrictions of $S_C$ and $S_B$ (or rather of their extensions to the corresponding multiplier algebras).

\bf \snl Proof\rm:
Denote $S_C(C)$ by $B'$ and $S_B(B)$ by $C'$. Then $S_B$ is an anti-isomorphism from $B$ to $B'$ and $S_C$ is an anti-isomorphism from $C$ to $C'$. By considering $E'$ as $\zeta(S_B\ot S_C)E$, we get in a straightforward way that $E'$ is a separability idempotent in $M(B'\ot C')$.
\hfill $\square$ \einspr

Remark that $E'$ is not necessarily a separability idempotent in the (possibly) larger algebra $M(B\ot C)$. In the regular case we study in the next section, we have $B'=B$ and $C'=C$ and therefore $E'$ is a separability in $M(B\ot C)$. An application of Proposition 1.14 will give $E'=E$, that is  $(S_B\ot S_C)E=\zeta E$, see Proposition 2.6 in the next section. In fact, we will see that this property is equivalent with $E$ being regular.
\snl
In the general case, it is a bit remarkable that on the one hand, $E$ is a (possibly) smaller idempotent than $E'$ (as $EE'=E'E=E$) while on the other hand, $E'$ is an element in a (possibly) smaller algebra (as $M(S_C(C)\ot S_B(B))\subseteq M(B\ot C)$). This is not yet very well understood. The problem is related with a similar, still open problem in the theory of possibly non-regular weak multiplier Hopf algebras. See a remark in Section 5 where we discuss conclusions and further research.

\newpage

\bf 2. Regular separability idempotents \rm
\nl
In this section, we will study regularity for separability idempotents. But before we start with the treatment of regular separability idempotents, we first consider an intermediate property.
\nl
\it Semi-regular separability idempotents \rm

\inspr{2.1} Definition \rm Let $B$ and $C$ be non-degenerate algebras and $E$ a separability idempotent in $M(B\ot C)$. Assume that not only $(b\ot 1)E$ and $E(1\ot c)$ belong to $B\ot C$ (as in the general case) but that also
$$E(b\ot 1)
\qquad\quad\text{and}\qquad\quad
(1\ot c)E$$
belong to $B\ot C$ for all $b\in B$ and $c\in C$. Then we call $E$ {\it semi-regular}.
\hfill $\square$ \einspr

We have the following characterization of this case.

\inspr{2.2} Proposition \rm
A separability idempotent $E$ in $M(B\ot C)$ is semi-regular if and only if the antipodal maps $S_B$ and $S_C$ have range in $C$ and $B$ respectively.

\bf \snl Proof\rm: 
i) Assume first that $E$ is semi-regular as in Definition 2.1. For all $b$ in $B$ we have $E(1\ot S_B(b))=E(b\ot 1)$ and this belongs to $B\ot C$ by assumption. If we now apply the functional $\varphi_B$, obtained in Proposition 1.11 on the first factor, we obtain 
$$S_B(b)=(\varphi_B\ot\iota)(E(b\ot 1))\in (\varphi_B\ot\iota)(B\ot C)\subseteq C.$$
Similarly $S_C(c)\in B$ for all $c\in C$.
\snl
ii) Conversely, let $b\in B$ and assume that $S_B(b)\in C$. Again from $E(b\ot 1)=E(1\ot S_B(b))$, we now get that $E(b\ot 1)\in B\ot C$. Similarly, if $c\in C$ and $S_C(c)\in B$, we  have $(1\ot c)E=(S_C(c)\ot 1)E\in B\ot C$. This proves that $E$ is semi-regular.
\hfill $\square$ \einspr

\it Regular separability idempotents \rm
\nl
We will see that there are various ways to characterize {\it regular} separability idempotents. We have chosen for the following definition.

\inspr{2.3} Definition \rm Let $B$ and $C$ be non-degenerate algebras and $E$ a separability idempotent in $M(B\ot C)$. We call $E$ {\it regular} if $\zeta E$ is a separability idempotent in $M(C\ot B)$. Recall that we use $\zeta$ for the flip map from $B\ot C$ to $C\ot B$, as well as for its unique extension to $M(B\ot C)$.
\hfill $\square$ \einspr

It is not hard to see that this is the same as assuming that $E$ is also a separability idempotent in $M(B^{\text{op}}\otimes C^{\text{op}})$.
\snl
It follows from the definition that in the regular case, all elements of the form 
$$\align  
&E(1\ot c) \qquad\qquad\text{and}\qquad\qquad (b\ot 1)E \\
&(1\ot c)E \qquad\qquad\text{and}\qquad\qquad E(b\ot 1) 
\endalign$$
belong to $B\ot C$ for all $b\in B$ and $C\in C$. This means that $E$ is semi-regular as in Definition 2.1. 
\snl
Moreover now there also exist non-degenerate anti-homomorphisms $S'_B:B\to M(C)$ and $S'_C:C\to M(B)$ satisfying and uniquely characterized by
$$(b\ot 1)E=(1\ot S'_B(b))E
\qquad\quad\text{and}\qquad\quad
E(1\ot c)=E(S'_C(c)\ot 1)$$
for all $b\in B$ and $c\in C$.
\snl
As regularity implies semi-regularity, it follows that the maps $S_B$ and $S'_B$ have range in $C$ and that the maps $S_C$ and $S'_C$ have range in $B$.
\snl
As an immediate consequence, we get the following characterization of a regular separability idempotent.

\inspr{2.4} Proposition \rm
A separability idempotent is regular if and only if the original anti-homomorphisms $S_B$ and $S_C$ have range in $C$ and $B$ respectively and are surjective. So they both are anti-isomorphisms.

\bf \snl Proof\rm:
i) First assume that $E$ is regular. Then for all $b\in B$ we have
$$E(b\ot 1)=E(1\ot S_B(b))=E(S'_C(S_B(b))\ot 1)$$
and it follows that $S'_C(S_B(b))=b$.  
So $S'_C$ is surjective and as we know that it is injective, it is a bijective map. Moreover, because  $S'_C(S_B(b))=b$ for all $b$, it follows that also $S_B$ is bijective from $B$ to $C$ and that $S_B$ and $S'_C$ are each others inverses.
\snl
A similar argument will give that also the two other maps $S_C$ and $S'_B$ are bijective and each others inverses.
\snl
ii) Conversely, assume that $S_B$ and $S_C$ are bijective from $B$ to $C$ and from $C$ to $B$ respectively. From $E(b\ot 1)=E(1\ot S_B(b))$ it follows that $E(b\ot 1)\in B\ot C$ for all $b$. Similarly, also $(1\ot c)E\in B\ot C$ for all $c\in C$. Further we simply define $S'_C$ as the inverse of $S_B$ and $S'_B$ as the inverse of $S_C$. Then we obtain that $\zeta E$ is again a separability idempotent in $M(C\ot B)$. 
\snl
This completes the proof.
\hfill $\square$ \einspr

Observe the difference with semi-regularity. In that case, we only have that $S_B$ has range in $C$ and $S_C$ has range in $B$, but it is not necessarily true that their ranges are all of $C$ and $B$ respectively as in the regular case.
\snl
There is also the following important remark.

\inspr{2.5} Remark \rm
In the earlier version of this paper (cf.\ [VD3.v1]), only regular separability idempotents were considered and for that reason, they were simply called separability idempotents. New results for the non-regular case, as defined and obtained in the previous section of this paper motivate this new point of view.
\snl
Moreover, in [VD-W4.v2] it is now shown that the canonical  idempotent of a (possibly non-regular) weak multiplier Hopf algebra is a separability idempotent in the sense of the previous section and that it is regular precisely if the original weak multiplier Hopf algebra is a regular one. Again, in an earlier version of that paper (see [VD-W4.v1]), only the result for regular weak multiplier Hopf algebras was obtained.
\snl 
We feel that this justifies the new terminology.
\hfill $\square$ \einspr

Recall that in general, we could prove that $EE'=E$ and $E'E=E$ where $E'=\zeta(S_B\ot S_C)E$ (see a remark made before Proposition 1.15 in Section 1). This means that $E\leq E'$. In the regular case, we have equality.

\inspr{2.6} Proposition \rm
For a regular separability idempotent, we have $\zeta(S_B\ot S_C)E=E$.

\snl\bf Proof\rm:
Define $E'=\zeta(S_B\ot S_C)E$ in $M(B\ot C)$ as before. We will use Proposition 1.14 to obtain that $E'=E$. We have to show that the antipodal maps for $E'$ are the same as for $E$. To do this, take $b\in B$. Then
$$\align E'(b\ot 1)&=(\zeta(S_B\ot S_C)E)(b\ot 1)\\
&=\zeta(S_B\ot S_C)((1\ot S_C^{-1}(b))E)\\
&=\zeta(S_B\ot S_C)((b\ot 1)E)\\
&=(\zeta(S_B\ot S_C)(E))(1\ot S_B(b))=E'(1\ot S_B(b)).
\endalign$$
Similarly we have that $S_C$ is the same for $E'$. So from Proposition 1.14 it follows that $E'=E$.
\hfill $\square$ \einspr

This formula is in accordance with the fact that $B$ and $C$ are the left and the right legs of $E$ and that $S_B:B\to C$ and $S_C:C\to B$. It also seems to suggest that conversely, if this equality holds, that is if $E'=E$ and not only $E\leq E'$, then $S_B$ and $S_C$ map $B$ and $C$ onto each other and so that $E$ is regular.
\snl
This is indeed the case as we show in the following proposition.

\inspr{2.7} Proposition \rm
If $E$ is a separability idempotent satisfying $\zeta(S_B\ot S_C)E=E$, then it is regular.

\snl\bf Proof\rm:
We will use the Sweedler type notation for $E$.
\snl
In Proposition 1.9 we have seen that $E_{(1)}\ot bS_C(E_{(2)})$ belongs to $B\ot B$ for all $b$ in $B$. 
\snl
Now assume that $E=S_C(E_{(2)})\ot S_B(E_{(1)})$. Then 
$$(b\ot 1)E=bS_C(E_{(2)})\ot S_B(E_{(1)})\tag"(2.1)"$$
and this belongs to $B\ot S_B(B)$ by the remark above. This implies that $C\subseteq S_B(B)$. 
\snl
On the other hand, assume that $\omega$ is a linear functional on $M(C)$ that vanishes on $C$. Using (2.1) we find $bS_C(E_{(2)})\omega(S_B(E_{(1)}))=0$ for all $b$. Replace again $b$ by $b'S_C(c)$. It follows now that $b'S_C(E_{(2)}c)\omega(S_B(E_{(1)}))=0$ for all $b'$ and all $c$. We can cancel $b'$ and use that $S_C$ is injective to obtain that also
$E_{(2)}c\omega(S_B(E_{(1)}))=0$ for all $c$. This implies that $\omega(S_B(b))=0$ for all $b$. Therefore $S_B(B)\subseteq C$.
\snl
As a consequence, $S_B(B)=C$. In a similar way we find $S_C(C)=B$. This means that $E$ is regular by Proposition 2.4 and the proof is complete.
\hfill $\square$ \einspr

We have now  different characterizations of regular separability idempotents. 
\nl
\it The distinguished linear functionals $\varphi_B$ and $\varphi_C$ in the regular case \rm
\nl
We now consider the linear functionals $\varphi_B$ and $\varphi_C$ as obtained in Proposition 1.11. 
\snl
From the equality $\zeta(S_B\ot S_C)E=E$ and the uniqueness of the linear functionals, it  follows that here $\varphi_C=\varphi_B\circ S_C$ and $\varphi_B=\varphi_C\circ S_B$. Consequently, $\varphi_C$ is invariant under the automorphism $S_BS_C$ of $C$ and $\varphi_B$ is invariant under the automorphism $S_CS_B$ of $B$. This last result also follows from the KMS-property in the next proposition.

\inspr{2.8} Proposition \rm
Assume again that $E$ is regular.
The linear functionals $\varphi_B$ and $\varphi_C$ are faithful. Moreover they satisfy
$$\varphi_B(bb')=\varphi_B(b'\sigma_B(b))
\qquad\quad\text{and}\qquad\quad
\varphi_C(cc')=\varphi_C(c'\sigma_C(c))$$
for all $b,b'\in B$ and $c,c'\in C$ where $\sigma_B=(S_CS_B)^{-1}$  and $\sigma_C=S_BS_C$.
\hfill $\square$ \einspr

The faithfulness can be obtained as in the proof of Proposition 1.12. It also follows from the one-sided faithfulness properties because $E$ is also a separability idempotent in \newline $M(B^{\text{op}}\ot C^{\text{op}})$. The KMS-properties have essentially been obtained already in Proposition 1.13. Now however, we know that the maps $S_CS_B$ and $S_BS_C$ are automorphisms of $B$ and $C$ respectively.
\snl
In the finite-dimensional case, for a faithful linear functional, the modular automorphism always exists (i.e.\ in our terminology, a faithful functional always has the KMS-property). In that case, the inverse of the modular automorphism is the Nakayama automorphism. In the general case, this is no longer true. The existence of the modular automorphisms is an extra property.
\snl
Recall also that the modular automorphisms are always trivial on central elements. So, if $z$ is central in $B$, it follows that $S_C(S_B(z))=z$. In fact, this can be seen directly because when $z$ is central, from $E(z\ot 1)=E(1\ot S_B(z))$ it follows that also $(z\ot 1)E=(1\ot S_B(z))E$ so that $S_C(S_B(z))=z$. 
We will use this later (in Theorem 4.13 of Section 4).
\nl
In the next section we consider some special cases. But we finish this section with a result, well-known in the finite-dimensional case.
\snl
Assume that $E$ is regular and let $F=(\iota\ot S_C)E$. It is a well-defined element in $M(B\ot B^{\text{op}})$ as we saw in Section 1. In Proposition 1.9 we have seen that $F(1\ot b)\in B\ot B$. Because $E$ is now assumed to be regular, we have that both $F(1\ot B)$ and $(B\ot 1)F$ are subsets of $B\ot B$. Moreover, in this case, 
$$(b\ot 1)F=F(1\ot b)\tag"(2.2)"$$ 
for all $b$ in $B$ because $S_C(C)=B$ (see Proposition 2.4). This yields the following result. 

\inspr{2.9} Proposition \rm
Assume that $E$ is regular. Define 
$\Delta:B\to B\ot B$ by 
$$\Delta(b)=(b\ot 1)F.$$ 
Then $\Delta$ is coassociative and the distinguished linear functional $\varphi_B$ is a counit.
\snl\bf Proof\rm:
i) Using the Sweedler notation $E_{(1)}\ot E_{(2)}$ for $E$ we have by definition $\Delta(b)=bE_{(1)})\ot S_C(E_{(2)})$. We find
$$\align
(\Delta\ot\iota)\Delta(b)&=(bE_{(1)}\ot 1)F)\ot S_C(E_{(2)}) \tag"(2.3)"\\
(\iota\ot\Delta)\Delta(b)&=(b\ot 1)F(1\ot E_{(1)})\ot S_C(E_{(2)})\tag"(2.4)"
\endalign$$
These expressions are the same, using (2.2). This proves coassociativity.
\snl
ii) If we apply $\varphi_B\ot \iota$ on $F(1\ot b)$ we find $b$ because  $(\varphi_B\ot \iota)E=1$. On the other hand, if we apply $\iota\ot\varphi_B$ to $(b\ot 1)F$, we also get $b$ now because $\varphi_B\circ S_C=\varphi_C$ (see a remark before Proposition 2.7). 
This proves that $\varphi_B$ is a counit on the coalgebra $(B,\Delta)$.
\hfill $\square$ \einspr

We clearly can also consider a coproduct $\Delta$ on $C$ defined by 
$\Delta(C)=((S_B\ot\iota)(E))(1\ot c)$. 
Now the counit is $\varphi_C$.
\snl
From this point of view, it would be meaningful to call the distinguished linear functionals $\varphi_B$ and $\varphi_C$ the counits. But then one has to keep in mind that they are counits for the coproducts defined as above.
\nl
In the non-regular case, we can still define these coproducts, but they will not map into the tensor product but rather into the multiplier algebras.
\snl
However, coassociativity only seems to be true in the regular case. Indeed, if e.g.\ we define $\Delta$ on $B$ as above by 
$\Delta(B)=(b\ot 1)F$ where again $F=(\iota\ot S_C)E$, we get $\Delta(b)\in B\ot M(B)$. In order to have coassociativity, we need equality of the expressions on the right in 
(2.3) and (2.4). Because the antipodal map $S_C$ is injective, this means that (using two copies of $E$ as we have done before)
$$
E_{(1)}E^\circ_{(1)}\ot S_C(E^\circ_{(2)})\ot E_{(2)}
=E^\circ_{(1)}\ot S_C(E^\circ_{(2)})E_{(1)}\ot E_{(2)}\tag"(2.5)"
$$
We rewrite the left hand side as 
$$
E_{(1)}\ot S_C(E^\circ_{(2)})\ot E_{(2)}S_B(E^\circ_{(1)}).
$$
If we now apply the functional $\varphi_B$ on the first factor, then the equality (2.5) will give 
$$
S_C(E^\circ_{(2)})\ot S_B(E^\circ_{(1)})=E_{(1)}\ot E_{(2)}.
$$
But this means $S_C(E_{(2)})\ot S_B(E_{(1)}))=E$ and we have seen before, in Proposition 2.6, that this implies regularity.
\snl
Remark that a similar argument will show that the equality $(b\ot 1)F=F(1\ot b)$ which is always true for $b\in S_C(C)$ will be true for all $b\in B$ if and only if $E$ is regular.
\nl
In the next section, we consider examples and special cases. And we will see that regularity is automatic in some of these special cases.
\nl\nl

\bf 3. Special cases and examples \rm 
\nl
We will now treat {\it some special cases} and give  {\it examples}. All of these are examples of regular separability idempotents. In Section 5, where we draw conclusions and discuss possible further research, we will explain why it seems difficult to find examples of non-regular separability idempotents. 
\nl
\it The finite-dimensional case \rm
\nl
In this item, we consider the case where the underlying algebras $B$ and $C$ are finite-dimensional.
\snl
Assume that we have a separability idempotent $E\in M(B\ot C)$ and that the algebras $B$, $C$ are finite-dimensional. Because they both have local units, they must be unital. Then $E\in B\ot C$. And regularity is automatic as we argue next.

\inspr{3.1} Proposition \rm
Assume that $B$ and $C$ are finite-dimensional unital algebras and that $E$ is a separability element in $B\ot C$. Then $S_B$ is an anti-isomorphism from $B$ to $C$ and $S_C$ is an anti-isomorphism from $C$ to $B$. Hence $E$ is regular.

\snl\bf Proof\rm:
Because the algebras are unital we have $S_B:B\to C$ and $S_C:C\to B$. Both maps are injective. This means that the dimension of $B$ is less or equal to the dimension of $C$ and vice versa. Therefore the dimensions are the same and the maps are bijections. The result now follows from the characterization of regularity that we found in Proposition 2.4 of the previous section.
\hfill $\square$ \einspr

In the introduction we gave already a simple example in the finite-dimensional setting. Now we make it somewhat more complicated.
\snl
First recall the example from the introduction. For $B$ and $C$ we have two copies of the algebra $M_n(\Bbb C)$ of $n\times n$ matrices over the complex numbers. We denote a set of matrix elements by $(e_{ij})$ where $i,j=1,2,\dots,n$. Then we have the idempotent
$$E_0=\frac1n \sum_{i,j=1}^n e_{ij}\ot e_{ij}$$
in $B\ot C$. It is a separability idempotent in the sense of this paper. The anti-isomorphisms $S_B$ and $S_C$ are the transposition of matrices, given by $S_0: e_{ij}\mapsto e_{ji}$ for all $i,j$. In this case $S_B$ and $S_C$ are each others inverses. The distinguished functionals $\varphi_B$ and $\varphi_C$ are nothing else but $n\text{Tr}$ where $\text{Tr}$ is the trace on $M_n(\Bbb C)$, normalized so that $\text{Tr}(1)=n$. The modular automorphisms are trivial. This is a consequence of the fact that the linear functionals are traces. And of course it is also seen from the formulas for the modular automorphisms in general because here $S_B$ and $S_C$ are each others inverses.
\snl
We now modify this example.

\inspr{3.2} Proposition \rm 
Let $n\in \Bbb N$ and take again for $B$ and $C$ the algebra $M_n(\Bbb C)$ as above. Let $r$ and $s$ be any two invertible matrices in $M_n(\Bbb C)$ and assume that $\text{Tr}(sr)=n$ where $\text{Tr}$ is the trace on $M_n(\Bbb C)$, normalized so that $\text{Tr}(1)=n$ as before. Put
$$E=(r\ot 1)E_0(s\ot 1)$$
where $E_0$ is as above. Then $E$ is a regular separability idempotent in $B\ot C$. The antipodal maps are given by
$$S_B(b)=S_0(sbs^{-1})
\qquad\quad\text{and}\qquad\quad
S_C(c)=rS_0(c)r^{-1}$$
for all $b$ in $B$ and $c$ in $C$ where $S_0$ is transposition of matrices. The distinguished linear functionals $\varphi_B$ and $\varphi_C$ are given by
$$\varphi_B=n\text{Tr}(p\,\cdot\,)
\quad\quad\text{and}\quad\quad
\varphi_C=n\text{Tr}(q\,\cdot\,)$$
where $p=(rs)^{-1}$ and $q=S_0(sr)^{-1}$. 

\bf \snl Proof\rm:
i) We find
$$E^2=(r\ot 1)E_0(sr\ot 1)E_0(s\ot 1)$$
and because
$$\align 
E_0(sr\ot 1)E_0
&=\frac{1}{n^2}\sum_{i,j,\ell,k} e_{ij}sre_{\ell k}\ot e_{ij}e_{\ell k}\\
&=\frac{1}{n^2}\sum_{i,j,k} e_{ij}sre_{jk}\ot e_{ik}\\
&=\frac{1}{n^2}\text{Tr}(sr)\sum_{i,k} e_{ik}\ot e_{ik}= E_0,
\endalign$$
we see that actually $E$ is still an idempotent in $M_n(\Bbb C)\ot M_n(\Bbb C)$. 
\snl
ii) A straightforward calculation gives the existence of the antipodal maps, as well as the formulas for $S_B$ and $S_C$. Indeed, if $b\in B$ we have
$$E(b\ot 1)=(r\ot 1)E_0((sbs^{-1})s\ot 1)=(r\ot 1)E_0(s\ot S_0(sbs^{-1}))=E(1\ot S_0(sbs^{-1})).$$
Similarly for $c\in C$ we find
$$(1\ot c)E=(r\ot c)E_0(s\ot 1)=(rS_0(c)\ot 1)E_0(s\ot 1)=(rS_0(c)r^{-1}\ot 1)E$$
\snl
iii) 
Let $\varphi_B$ be given by $\varphi_B=n\text{Tr}(p\,\cdot\,)$ as in the formulation of the result.
We have
$$\align
(\text{Tr}\ot\iota)((p\ot 1)E))
&=(\text{Tr}\ot\iota)((pr\ot 1)E_0(s\ot 1))\\
&=(\text{Tr}\ot\iota)((s^{-1}\ot 1)E_0(s\ot 1))\\
&=(\text{Tr}\ot\iota)E_0\\
&=\frac1n\sum_{ij}\text{Tr}(e_{ij})e_{ij}
=\frac1n\sum_j e_{jj}.
\endalign$$
Therefore $(\varphi_B\ot\iota)E=1$. Similarly we find $(\iota\ot\varphi_C)E=1$ when $\varphi_C$ is defined as $n\text{Tr}(q\,\cdot\,)$.
\snl
iv)
The modular automorphisms $\sigma_B$ and $\sigma_C$ are given by
$$\sigma_B(b)=pbp^{-1}
\qquad\quad\text{and}\qquad\quad
\sigma_C(c)=qcq^{-1}$$
where $p$ and $q$ are defined as above. This is an immediate consequence of the expressions obtained for $\varphi_B$ and $\varphi_C$. The formulas can also be obtained from the general formulas for $\sigma_B$ and $\sigma_C$ (obtained in Proposition 2.8), given the formulas for $S_B$ and $S_C$.
\hfill$\square$\einspr

If we allow automorphisms of the factors, we can even assume in this case that either $r$ or $s$ is equal to $1$. 
\snl
Remark that in the case where $\text{Tr}(sr)=0$, so that $E^2=0$, we still have the anti-isomorphisms $S_B$ and $S_C$, as well as the functionals $\varphi_B$ and $\varphi_C$ and the associated modular automorphisms. This suggests that there might be a more general and still interesting theory where full elements $E$ in $M(B\ot C)$ are considered, without the assumption of having an idempotent element. We refer to a similar remark (after Proposition 3.12) when we consider the case of a discrete quantum group later in this section. See also Section 5 where we discuss possible future research.
\snl 
This example is very important as we see from the following result.

\inspr{3.3} Proposition \rm
Any separability idempotent in $M_n(\Bbb C)\ot M_n(\Bbb C)$ is as in the previous example.

\snl \bf Proof\rm:
Assume that $E$ is a separability idempotent in $M_n(\Bbb C)\ot M_n(\Bbb C)$. We know that it is regular and therefore the antipodal maps are anti-isomorphisms of $M_n(\Bbb C)$. It follows that here exist invertible matrices $s$ and $r$ such that
$$S_C(c)=rS_0(c)r^{-1}
\qquad\quad\text{and}\qquad\quad
S_B(b)=S_0(sbs^{-1})$$
for all $b,c$. Here again, $S_0$ is transposition of matrices. Then $E$ has the same antipodal maps as the separability idempotent  associated with these invertible matrices as in Proposition 3.2 above. By Proposition 1.14, they are the same. Because $E^2=E$, we need to have that $\text{Tr}(rs)\neq 0$. Then we can, if necessary, scale $r$ or $s$ to arrive at the condition $\text{Tr}(rs)=n$. This completes the proof.
\hfill$\square$\einspr

The example above is also important because any regular separability idempotent is a direct sum of a possibly infinite numbers of separability idempotents of the form above (see Theorem 4.13 in Section 4).
\nl
\it The abelian case \rm
\nl
Let us first make the following observation.

\inspr{3.4} Lemma \rm
If the algebra $B$ is abelian, then also $C$ has to be abelian. And vice versa.

\snl\bf Proof\rm: 
Because $B$ is abelian, also $M(B)$ is abelian and for all $c,c'$ we will have
$$S_C(cc')=S_C(c')S_C(c)=S_C(c)S_C(c')=S_C(c'c).$$
As $S_C$ is injective it follows that $cc'=c'c$. Hence $C$ is also abelian.
\hfill $\square$ \einspr

\inspr{3.5} Proposition \rm
Assume that $B$ and $C$ are abelian non-degenerate algebras. If $E$ is a separability idempotent in $M(B\ot C)$, then it is regular. Moreover $S_B$ and $S_C$ are each others inverses.

\bf \snl Proof\rm:
First remark that $E$ is semi-regular. Furthermore, for $b$ in $B$ we have
$$\align
E(b\ot 1)&=E(1\ot S_B(b))=(1\ot S_B(b))E\\
&=(S_C(S_B(b))\ot 1)E=E(S_C(S_B(b))\ot 1).\endalign$$
Then $S_C(S_B(b))=b$ and as in the proof of Proposition 2.4 it follows that $S_B$ and $S_C$ are bijections from $B$ to $C$ and from $C$ to $B$ respectively and that they are each others inverses.
\hfill $\square$ \einspr

We could also use the characterization that $E$ is regular if and only if $E$ is also a separability idempotent in $M(B^{\text{op}}\ot C^{\text{op}})$ because now $B^{\text{op}}=B$ and $C^{\text{op}}=C$.
\nl
Consider here the following trivial but basic example.

\inspr{3.6} Example \rm
Let $X$ be any set. Let $B$ and $C$ both be the algebra $K(X)$ of complex functions with finite support in $X$ and pointwise operations. This is a non-degenerate, abelian algebra. It has an identity if and only if $X$ is a finite set. In general, it clearly has local units. 
The algebra $B\ot C$ is $K(X\times X)$, the algebra of complex functions with finite support on the Cartesian product $X\times X$ of $X$ with itself. The multiplier algebra $M(B\ot C)$ is identified with the algebra $C(X\times X)$ of all complex functions on this Cartesian product.
\snl
Define $E\in M(B\ot C)$ by $E(x,y)=1$ if $x=y$ and $E(x,y)=0$ otherwise for $x,y\in X$. Then $E$ is a regular separability idempotent and the antipodal maps are the identity maps.
\snl
Indeed, first of all, it is clear that $E$ is idempotent. Moreover, if $f$ is any function on $X$ with finite support, then also $E(1\ot f)$ will have finite support in $X\times X$. The same is true for the three other expressions of this type. If $f$ is the function $\delta_x$ that is one in a given point $x$ and $0$ in all other points, then 
$E(1\ot f)=\delta_x\ot \delta_x$. Similarly for $E(f\ot 1)$. This shows that $E$ is full. Finally, for all $f\in K(X)$ we trivially have $E(f\ot 1)=E(1\ot f)$ proving that the antipodal maps exist and are the identity maps.
\snl
The linear functionals $\varphi_B$ and $\varphi_C$ coincide and are given by $f\mapsto \sum_x f(x)$. The modular automorphisms are of course trivial.
\hfill $\square$ \einspr

Later, we will see that any separability idempotent in the abelian case has to be of this form. This is a consequence of the structure theorem, Theorem 4.13 in Section 4.
\nl
\it The involutive case \rm
\nl
In this item, we assume that the underlying algebras $B$ and $C$ are $^*$-algebras and that the separability idempotent $E$ in $M(B\ot C)$ is self-adjoint. Again, in this case regularity is automatic. We have also some other extra properties.

\inspr{3.7} Proposition \rm
Assume that $B$ and $C$ are $^*$-algebras and that $E$ is self-adjoint. Then $E$ is regular. Moreover, 
$$S_C(S_B(b)^*)^*=b
\qquad\quad\text{and}\qquad\quad
S_B(S_C(c)^*)^*=c$$ 
for all $b\in B$ and $c\in C$.

\snl\bf Proof\rm: 
By taking adjoints, we find that not only $E(1\ot c)$ and $(b\ot 1)E$ are in $B\ot C$, but also 
$(1\ot c)E$ and $E(b\ot 1)$ for all $b\in B$ and $c\in C$. Also, if we take adjoints of the equations
$$E(b\ot 1)=E(1\ot S_B(b))
\qquad\quad\text{and}\qquad\quad
(1\ot c)E=(S_C(c)\ot 1)E,$$
we find 
$$(b^*\ot 1)E=(1\ot S_B(b)^*)E
\qquad\quad\text{and}\qquad\quad
E(1\ot c^*)=E(S_C(c)^*\ot 1)$$
for all $b,c$. This proves the existence of the anti-homomorphisms $S'_B$ and $S'_C$ (as in Proposition 2.4) and that 
$$S'_B(b^*)=S_B(b)^*
\qquad\quad\text{and}\qquad\quad
S'_C(c^*)=S_C(c)^*$$
for all $b,c$. This proves that $E$ is regular. 
\snl
Because we also know that $S'_B$ and $S'_C$ are the inverses of $S_C$ and $S_B$ respectively, we find the formulas in the formulation of the proposition. This completes the proof.
\hfill $\square$ \einspr

In the involutive case, when $E$ is self-adjoint, it turns out that the distinguished linear functionals $\varphi_B$ and $\varphi_C$ are positive linear functionals. We expect that there is a direct proof of this result, based on the formula $E=E^*E$, but we have not been able to find this. We will need the structure theorem of Section 4, together with the results below.
\snl
Here is the prototype of an example in the involutive case. It is a special case of the example in Proposition 3.2.

\inspr{3.8} Example \rm
Let $n\in \Bbb N$ and take for $B$  and $C$ the $^*$-algebra $M_n(\Bbb C)$. Let $r$ be any  invertible element in $M_n(\Bbb C)$ and assume that $\text{Tr}(r^*r)=n$ where $\text{Tr}$ is the trace on $M_n(\Bbb C)$, normalized so that $\text{Tr}(1)=n$ as before. Now define
$$E=(r\ot 1)E_0(r^*\ot 1)$$
where 
$$E_0=\frac{1}{n}\sum_{i,j}e_{ij}\ot e_{ij}$$
as in the introduction. Again $(e_{ij})$ are matrix units in $M_n(\Bbb C)$. 
\snl
i) We find, just as in the proof of Proposition 3.2, 
$$E^2=(r\ot 1)E_0(r^*r\ot 1)E_0(r^*\ot 1)=(r\ot 1)E_0(r^*\ot 1)$$
because $Tr(r^*r)=n$. We see that $E$ is a self-adjoint idempotent in $M_n(\Bbb C)\ot M_n(\Bbb C)$. 
\snl
ii) As in Proposition 3.2, it is again straightforward to verify that $E$ is a separability idempotent in $M_n(\Bbb C)\ot M_n(\Bbb C)$. In fact, we can apply the result of this proposition. In this case, we get for the anti-isomorphisms $S_B$ and $S_C$ the formulas 
$$S(b)=S_0(r^*b{r^*}^{-1})
\tussenen
S_C(c)=rS_0(c)r^{-1}$$
for all $b$ and $c$. 
\snl
iii) The distinguished linear functionals are given by 
$$\varphi_B=n\text{Tr}(p\,\cdot\,)
\tussenen
\varphi_C=n\text{Tr}(q\,\cdot\,)$$
where $p=(rr^*)^{-1}$ and $q=S_0(r^*r)^{-1}$ . 
\snl 
iv) For the {\it modular automorphisms} we get
$$\sigma_B(b)=pbp^{-1}
\tussenen
\sigma_C(c)=qcq^{-1}$$
where $p=(rr^*)^{-1}$ and $q=S_0(r^*r)^{-1}$. 
\snl
v) The distinguished linear functionals are clearly positive as the elements $p$ and $q$ are positive.
\hfill$\square$\einspr

Also this example is universal in the following sense.

\inspr{3.9} Proposition \rm
Any self-adjoint separability idempotent in $M_n(\Bbb C)\ot M_n(\Bbb C)$ is  as in the previous example.

\snl\bf Proof\rm:
Let $E$ be any self-adjoint separability idempotent in $M_n(\Bbb C)\ot M_n(\Bbb C)$. Let $S_B$ and $S_C$ be the associated anti-isomorphisms. The anti-isomorphism $S_C$ has the form $S_C(c)=rS_0(c)r^{-1}$ for all $c$ for some invertible matrix $r\in M_n(\Bbb C)$. We know from the general theory that $S_B(b)=(S_C^{-1}(b^*)^*$ and then $S_B(b)=r^*S_0(b)r^{*-1}$ for all $b$. As in Proposition 3.3 we find that $E=(r\ot 1)E_0(r^*\ot 1)$. Finally, we can scale $r$ so that $\text{Tr}(r^*r)=n$ and this will complete the proof.
\hfill$\square$\einspr

We need to make the following remark. Consider again Example 3.8. Let $r=v|r|$ be the {\it polar decomposition} of $r$. So $|r|=(r^*r)^\frac12$, it is a positive self-adjoint invertible matrix and $v$ is a unitary matrix. Then we have $E=(\alpha_v\ot\iota)E'$ where $E'=(|r|\ot \iota)E_0(|r|\ot 1)$ and where $\alpha_v$ is the $^*$-automorphism $b\mapsto vbv^*$. It follows from this observation that we can assume in the previous proposition, up to a $^*$-isomorphism, that $r$ is a positive self-adjoint invertible matrix. This simplifies some formulas. We then have e.g.\ that $q=S_0(p)$.
\nl
\it Examples with direct sums \rm
\nl
We will now use the above examples to build new ones. First we consider a general construction.
\snl
Let $J$ be an index set. Assume that for all $\alpha\in J$ we have a pair of non-degenerate algebras $B_\alpha$ and $C_\alpha$. We consider the direct sums 
$$B=\sum_{\alpha\in J} B_\alpha
\qquad\quad\text{and}\qquad\quad
C=\sum_{\alpha\in J} C_\alpha.$$
Elements $b$ in $B$ are of the form $(b_\alpha)_{\alpha\in J}$ where $b_\alpha\in B_\alpha$ for all $\alpha\in J$ and where only a finite number of the components are non-zero. Similarly for $C$. There are natural identifications 
$$M(B)=\prod_{\alpha\in J}M(B_\alpha)
\qquad\quad\text{and}\qquad\quad
M(C)=\prod_{\alpha\in J}M(C_\alpha).$$
Elements in $M(B)$ have the form $(x_\alpha)_{\alpha\in J}$ where $x_\alpha\in M(B_\alpha)$ for all $\alpha\in J$ (without further restrictions) and similarly for elements in $M(C)$. 
For the tensor product $B\ot C$ and its multiplier algebra $M(B\ot C)$ we get 
$$\align B\ot C &= \sum_{\alpha,\beta}B_\alpha \ot C_\beta \\
M(B\ot C)&=\prod_{\alpha,\beta}M(B_\alpha\ot C_\beta).
\endalign$$
The basic result using these building blocks is now the following. We use the above notations.

\iinspr{3.10} Proposition \rm For all $\alpha\in J$ let $E_\alpha$ be a separability idempotent in $M(B_\alpha\ot C_\alpha)$. Define the element $E$ in $M(B\ot C)$ with components $E_\alpha$. Then $E$ is a separability idempotent in $M(B\ot C)$. The anti-homomorphisms $S_B$ and $S_C$ are given by 
$$S_B(b)_\alpha=S_B^\alpha(b_\alpha)
\tussenen
S_C(c)_\alpha=S_C^\alpha(c_\alpha)$$
 for all $b\in B$ and $c\in C$ (where of course $S_B^\alpha$ and $S_C^\alpha$ are the antipodal maps associated with the components $E_\alpha)$. The left and right functionals $\varphi_C$ on $C$ and $\varphi_B$ on $B$ are given by 
$$\varphi_B(b)=\sum_{\alpha \in J} \varphi_B^\alpha(b_\alpha)
\qquad\quad\text{and}\qquad\quad
\varphi_C(c) =\sum_{\alpha \in J} \varphi_C^\alpha(c_\alpha)$$
where $\varphi_B^\alpha$ and $\varphi_C^\alpha$ are the distinguished linear functionals for $E_\alpha$. Observe that these sums are finite sums as only finitely many terms are non-zero. The element $E$ is regular if and only if all components are regular.
\hfill$\square$\einspr

The proof is  straightforward and we leave it mostly as an exercise for the reader. We just make a few remarks. 
\snl
Clearly $E$ is an idempotent as all its components are idempotents. If $c\in C$, we find that $(E(1\ot c))_{\alpha,\beta}=0$ if $\alpha\neq\beta$ while $(E(1\ot c))_{\alpha,\alpha}=E_\alpha(1\ot c_\alpha)$ and this will be $0$ for all but finitely many $\alpha$. It follows that $E(1\ot c)\in B\ot C$. Similarly $(b\ot 1)E\in B\ot C$. To prove that $E$ is full, take e.g. a fixed index $\alpha_0$ and assume that $c$ is the element with $0$ in all components, except for the $\alpha_0$ component where it is an element $c_0$ of $C_{\alpha_0}$. For any functional $\omega$ on $C$ we will find that $(\iota\ot\omega)(E(1\ot c))$ will only have one non-zero component in the index $\alpha_0$, given by $(\iota\ot\omega_{\alpha_0})(E_{\alpha_0}(1\ot c_0))$ where $\omega_{\alpha_0}$ is the restriction of $\omega$ to the component $C_{\alpha_0}$. As we can reach all elements in $B_{\alpha_0}$ in such a way, we see that the left leg of $E$ is indeed all of $B$. Similarly for the right leg. 
\snl
It is obvious that $E(b\ot 1)=E(1\ot S_B(b))$ if we define $S_B$ on $B$ by $S_B(b)_\alpha=S_B^\alpha(b_\alpha)$ for all $\alpha$. Similarly for $S_C$ on $C$. This will show that we have a separability idempotent. The formulas for the functionals $\varphi_B$ and $\varphi_C$ are obvious. 
\snl
Finally, the modular automorphisms $\sigma_B$ and $\sigma_C$ are given by $\sigma_B(b)_\alpha=\sigma_B^\alpha(b_\alpha)$ and $\sigma_C(c)_\alpha=\sigma_C^\alpha(c_\alpha)$ for all $b\in B$ and $c\in C$ (where $\sigma_B^\alpha$ and $\sigma_C^\alpha$ are the modular automorphisms associated with the components $E_\alpha$). 
\snl
It is also clear that this construction will work in the involutive case as well. Indeed, if components are $^*$-algebras, then so will be $B$ and $C$ and if each $E_\alpha$ is self-adjoint, the same will be true for $E$.
\snl
This result can now be combined with the examples given in Proposition 3.2 and Example 3.8 to obtain infinite-dimensional cases. Indeed, we can apply Proposition 3.10 with the examples in Example 3.8 (in the involutive case) or with the example in Proposition 3.2 (in the general case). We will will see in the next section that this will provide all regular separability idempotents. In fact, the example in Example 3.6 is also of this type. We come back to this in Section 4, after the proof of Theorem 4.13.
\nl
\it Discrete quantum groups \rm
\nl
For our next example, we have the following. We refer to [VD1] and [VD2] for the theory of multiplier Hopf algebras and in particular to [VD-Z] for multiplier Hopf algebras of discrete type.

\iinspr{3.11} Proposition \rm
Let $(A, \Delta)$ be a regular multiplier Hopf algebra of discrete type.  Assume that the (left) cointegral $h$ is an idempotent. Then $\Delta(h)$ is a separability idempotent in $M(A\ot A)$. The anti-isomorphisms $S_B$ and $S_C$ coincide with the antipode of $A$. The distinguished linear functionals $\varphi_B$ and $\varphi_C$ are the right and left integrals on $A$ respectively.

\snl \bf Proof\rm: 
Recall that by definition, a left cointegral $h$ exists. However, it can happen that $h^2=0$ (when $\varepsilon(h)=0$ where $\varepsilon$ is the counit). If this is not the case, we can scale $h$ so that $\varepsilon(h)=1$. Then $h^2=h$ and $h$ will also be a right cointegral. This is what we assume here. We will use the term discrete quantum group for such multiplier Hopf algebras of discrete type with a normalized cointegral.
\snl
It follows that $\Delta(h)$ is an idempotent in $M(A\ot A)$. We know from the theory of discrete quantum groups that the legs of $\Delta(h)$ are all of $A$. This implies fullness of the idempotent.
\snl
For discrete quantum groups (where the left cointegral is also a right cointegral), it also holds
$$\Delta(h)(a\ot 1)=\Delta(h)(1\ot S(a))
\qquad\quad\text{and}\quad\qquad
(1\ot a)\Delta(h)=(S(a)\ot 1)\Delta(h)$$
for all $a$ where $S$ is the antipode of $A$. The first equality is true for any left cointegral and the second for any right cointegral. In this case, these are the same.
\snl
Because we assume that $(A,\Delta)$ is regular, the antipode $S$ is an anti-isomorphism of $A$ and it follows that we have a regular separability idempotent. It is clear from the formulas above that the antipode $S$ of $A$ gives the anti-isomor\-phisms $S_B$ and $S_C$.
\snl
And it is also clear the left and right integrals on $A$ give the distinguished linear functionals for the separability idempotent $\Delta(h)$.
\hfill$\square$\einspr

If $(A,\Delta)$ is a multiplier Hopf $^*$-algebra of discrete type, it is automatically regular, left and right cointegrals coincide and are self-adjoint idempotents. In this case, we get a self-adjoint separability idempotent. 
\snl
For a general multiplier Hopf algebra of discrete type $(A,\Delta)$ there is the possibility that the left and right cointegrals are different. In this case $h^2=0$ and $\Delta(h)$ is not idempotent. Still, some properties of a separability idempotent remain true and again, this suggest that a more general theory is feasible, just like we remarked before after we had the example in Proposition 3.2. See some remarks about this in Section 5.
\nl
\it Weak multiplier Hopf algebras \rm
\nl
Now we consider the case of weak multiplier Hopf algebras. We refer to [VD-W2], [VD-W3] and [VD-W4.v2] for the theory of weak multiplier Hopf algebras. We have the following result.

\iinspr{3.12} Proposition \rm
Let $(A,\Delta)$ be a  weak multiplier Hopf algebra. Denote by $B$ the image $\varepsilon_s(A)$ of $A$ under the source map $\varepsilon_s$ and by $C$ the image $\varepsilon_t(A)$ of $A$ under the target map $\varepsilon_t$. Then the canonical idempotent $E$ is a separability idempotent in $M(B\ot C)$. The antipode $S$ of $A$, when restricted to $B$ and $C$ gives the associated anti-isomorphisms $S_B$ and $S_C$ respectively. The left and right distinguished linear functionals on $C$ and $B$ are given by the maps
$$
\varphi_B:\varepsilon_s(a)\mapsto \varepsilon(a)
\qquad\quad\text{and}\quad\qquad
\varphi_C\varepsilon_t(a)\mapsto \varepsilon(a)$$ 
here $\varepsilon$ is the counit of $A$. The separability idempotent is regular if the original weak multiplier Hopf algebra is regular. The modular automorphisms are the restrictions of the square $S^2$ of the antipode.
\hfill\einspr

This can be found in  [VD-W4.v2]. Indeed, in Theorem 2.18 of [VD-W4.v2], it is shown that $E$ is a separability idempotent in $M(B\ot C)$ where $B$ is the image of the source map $\varepsilon_s$ and $C$ the image of the target map $\varepsilon_t$.The antipodal maps are obtained by first extending the antipode to all of $M(A)$ and then as the restrictions of this extension to $B$ and $C$ respectively. Finally, in Proposition 2.20 of [VD-W4.v2], the formulas for the distinguished linear functionals $\varphi_B$ and $\varphi_C$ are obtained. The statement about regularity is found in Proposition 2.19 of [VD-W4.v2].
\snl
Remark that the functionals $\varphi_B$ on $B$ and $\varphi_C$ on $C$ also exist when there are no integrals on $A$. This is different from the situation we had in the previous example with a discrete quantum group. In that case, the integrals on the legs of $E$ where the integrals on the original algebra $A$.
\snl
Also remark that this example is in fact essentially the general situation. Indeed, if we  have a separability idempotent $E\in M(B\ot C)$, we can make the algebra $C\ot B$ into a regular weak multiplier Hopf algebra so that the canonical idempotent is $1\ot E\ot 1$ where the first $1$ is the identity in $M(C)$ and the second $1$ the identity in $M(B)$. For this example, one has precisely that $1\ot B$ and $C\ot 1$ are the images of the source and target maps respectively. See Proposition 3.2 in [VD-W4.v2].  
\nl\nl

\bf 4. Modules and the structure of the base algebras \rm
\nl
In this section we study modules over the base algebras $B$ and $C$. We use the existence of the separability idempotent to obtain the complements of submodules. This eventually will yield results about the structure of the base algebras.
\snl
In the first place, we treat the general case but nicer results are obtained if the separability idempotent is regular. In fact, in this case, the results are very similar to the well-known results in the finite-dimensional case and more generally, for unital algebras. 
\snl
We start with recalling some notions and results about modules over non-degenerate, not necessarily unital algebras.
\nl
\it Unital and non-degenerate modules \rm
\nl
First recall the following definition (see Definition 3.1 and Proposition 3.2 in [Dr-VD-Z]).

\inspr{4.1} Definition \rm
Let $A$ be a non-degenerate algebra and $X$ a left $A$-module. We write the action as $(a,x)\mapsto ax$ for $a\in A$ and $x\in X$. Then $X$ is called a {\it unital module} if any element in $X$ is a sum of elements of the form $ax$ where $a\in A$ and $x\in X$. This condition is written as $X=AX$. The module is called {\it non-degenerate} if given $x\in X$, we have $x=0$ if $ax=0$ for all $a\in A$.
\hfill $\square$ \einspr

If the algebra $A$ has left local units and if $X$ is unital, then for all $x\in X$ there is an element $e\in A$ so that $ex=x$. This implies that the module is automatically non-degenerate.
\snl
There are similar notions and results for right $A$-modules and $A$-bimodules. A bimodule is called unital if it is unital both as a left and as a right module. Similarly, it is called non-degenerate if that is true when considered as a left module and as a right module. If the algebra has two-sided local units, a unital bimodule will automatically be non-degenerate. This is the case for the algebras $B$ and $C$ that we will consider.
\snl
We will need the extension of the action to the multiplier algebra (see Proposition 3.3 in [Dr-VD-Z]).

\inspr{4.2} Proposition \rm
Assume that $X$ is a unital and non-degenerate left $A$-module. Then the action of $A$ has a unique extension to a unital action of $M(A)$.
\hfill $\square$ \einspr

The proof is based on the following two observations. 
\snl
If $x=\sum_i a_ix_i$ and $m\in M(A)$, then we must have $mx=\sum_i (ma_i)x_i$. This is used to prove uniqueness. It also suggest how to obtain the extension. Observe that 
$$b(\textstyle\sum_i (ma_i)x_i)=\textstyle\sum_i (bma_i)x_i=(bm)(\textstyle\sum_i a_ix_i)$$
for all $b\in A$ and because $X$ is assumed to be non-degenerate, we get $\sum_i (ma_i)x_i=0$ if $\sum_i a_ix_i=0$. This is used to prove the existence. It is also obvious that $mx=x$ when $m=1$.
\snl
The procedure can also be applied in the case of a unital, non-degenerate right module as well as to unital bimodules. 
\nl
\it Projection maps on submodules \rm
\nl
Let $B$ and $C$ be non-degenerate algebras and $E$ a separability idempotent in $M(B\ot C)$. In Proposition 1.6 we have obtained the formulas
$$bE_{(1)}S_C(E_{(2)})=b
\qquad\quad\text{and}\qquad\quad
S_B(E_{(1)})E_{(2)}c=c$$
for $b\in B$ and $c\in C$ where we use the Sweedler type notation $E_{(1)}\ot E_{(2)}$ for $E$ as before. Recall also that $B$  and  $C$ have local units (cf.\ Proposition 1.10).
\snl
For next proposition, we consider a {\it unital left $C$-module $X$} and a {\it submodule $Y$ of $X$}. Because $C$ has local units, the module $X$ is also non-degenerate.  For all elements $x$ in $X$ there is an element $e\in C$ so that $ex=x$. If we apply this to an element of $Y$, we see that the submodule is still unital and non-degenerate. We can extend the action of $C$ on $X$ to $M(C)$ as in Proposition 4.2. It is easily seen that $Y$ is still invariant under this action of $M(C)$. In fact, the action of $M(C)$ on $Y$ obtained in this way is  the same as the one we would obtain by applying the extension procedure for the action of $C$ on $Y$.

\inspr{4.3} Proposition \rm
Let $q:X\to Y$ be a linear projection map. Then we can define a new linear map $p:X\to Y$ by
$$p(x)=S_B(E_{(1)})q(E_{(2)}x)$$
for $x\in X$. Also $p$ is a projection map and now
$$p(S_B(b)x)=S_B(b)p(x)$$
for all $x\in X$ and $b\in B$.

\snl\bf Proof\rm:
i) We first show that $p$ is well-defined from $X$ to $Y$. 
\snl
Let $x$ be in $X$. Because $C$ has local units, we have an element $e$ in $C$ so that $ex=x$. Then
$$S_B(E_{(1)})\ot E_{(2)}ex\in S_B(B)\ot CX\subseteq M(C)\ot X$$ 
so that 
$$S_B(E_{(1)})q(E_{(2)}ex)$$
is well-defined as an element in $Y$. Now suppose that $e$ and $e'$ are elements in $C$ so that $ex=x$ and $e'x=x$. Choose $f$ in $C$ so that $fe=e$ and $fe'=e'$. Then
$$\align S_B(E_{(1)})q(E_{(2)}ex)
&=S_B(E_{(1)})q(E_{(2)}fex)\\
&=S_B(E_{(1)})q(E_{(2)}fe'x)\\
&=S_B(E_{(1)})q(E_{(2)}e'x).
\endalign$$
So we can define $p(x)=S_B(E_{(1)})q(E_{(2)}ex)$ when $e$ is an element in $C$ so that $ex=x$. We can safely write $p(x)=S_B(E_{(1)})q(E_{(2)}x)$.
\snl
ii) To show that $p$ is a linear map, we use that for two elements $x,x'$ in $X$, we can find a single $e$ in $C$ so that $ex=x$ and $ex'=x'$.
\snl
iii) We already know that $p$ maps $X$ to $Y$. On the other hand, if $x\in Y$, then 
$$S_B(E_{(1)})\ot E_{(2)}x\in S_B(B)\ot CY\subseteq M(C)\ot Y$$
and so, as $q$ is a projection map on $Y$,  
$$p(x)=S_B(E_{(1)})E_{(2)}x=x.$$
It follows that $p$ is also a projection map onto $Y$.
\snl
iv) Finally, let $b\in B$ and $x\in X$. Then
$$\align p(S_B(b)x)
&=S_B(E_{(1)})q(E_{(2)}S_B(b)x)\\
&=S_B(E_{(1)}b)q(E_{(2)}x)\\
&=S_B(b)S_B(E_{(1)})q(E_{(2)}x)\\
&=S_B(b)p(x).
\endalign$$
\vskip -0.8 cm
\hfill $\square$ \einspr

In the regular case, we have $S_B(B)=C$ and we find that $p(cx)=cp(x)$ for all $x\in X$ and all $c\in C$. It is not clear if this property is still true in the non-regular case, or even in the semi-regular case. Compare with the discussion at the end of Section 2 where we had a similar problem. There, the result was only true in the regular case.
\snl
In the semi-regular case, we have $S_B(B)\subseteq C$ and we find $p(cx)=cp(x)$ for all $x\in X$ but only for $c$ in the subalgebra $S_B(B)$ of $C$.
\nl
We have a completely similar result for right $B$-modules.

\inspr{4.4} Proposition \rm
Let $X$ be a unital right $B$-module and $Y$ a submodule. Let $q:X\to Y$ be a linear projection map. We can define a new map $p:X\to Y$ by
$$p(x)=q(xE_{(1)})S_C(E_{(2)})$$
for $x\in X$. Then $p$ is a projection map so that moreover
$$p(xS_C(c))=p(x)S_C(c)$$
for all $x\in X$ and $c\in C$.
\hfill $\square$ \einspr

Also here, in the regular case where $S_C(C)=B$ we will find $p(xb)=p(x)$ for all $x\in X$ and all $b\in B$. Again, it is not clear if this property is still true in the non-regular case. 
\nl
Next we consider the two other cases, a submodule of a right $C$-module and a submodule of a left $B$-module. Here the situation is a little more complicated. We will need that $((\iota\ot S_C)E)(1\ot b)$ belong to $B\ot B$ for all $b\in B$ and that $(c\ot 1)(S_B\ot\iota)E$ belong to $C\ot C$ for all $c\in C$ (cf.\ Proposition 1.9).

\inspr{4.5} Proposition \rm
Assume that $E$ is a separability idempotent in $M(B\ot C)$. Let $X$ be a unital right $C$-module and $Y$ a submodule. Assume that $q$ is a linear projection map from $X$ onto $Y$. Then we can define a projection map $p:X\to Y$ by
$$p(x)=q(xS_B(E_{(1)}))E_{(2)}$$
where $x\in X$. Moreover $p(xS_B(b))=p(x)S_B(b)$ for all $x\in X$ and $b\in B$.

\snl\bf Proof\rm:
The proof is very much the same as the one of Proposition 4.2. The main difference lies in showing that $p$ also here is well-defined. 
\snl
Because $X$ is assumed to be unital and because $C$ has local units, given $x\in X$ there is an element $e\in C$ so that $xe=x$. Then we can define
$$ q(x(eS_B(E_{(1)})))E_{(2)}$$
as an element in $Y$ because of the result in Proposition 1.9 we just mentioned. Again, this expression will not be dependent on the choice of $e$. Therefore we can define $p$.
The rest of the argument is also very similar as in the proof of Proposition 4.2.
\hfill $\square$\einspr

We have a similar result for unital left $B$-modules.

\inspr{4.6} Proposition \rm
Let $X$ be a unital left $B$-module and $Y$ a submodule. Assume that $q$ is a linear projection map from $X$ onto $Y$. Then we can define a projection map $p:X\to Y$ by
$$p(x)=E_{(1)})q(S_C(E_{(2)})x).$$
Moreover $p(S_C(c)x)=S_C(c)p(x)$ for all $x\in X$ and $c\in C$.
\hfill $\square$ \einspr

If we combine the two techniques, we get the following results for bimodules.

\inspr{4.7} Proposition \rm
Assume that $X$ is a unital $C$-bimodule and $Y$ a submodule. Then there is a projection map $p$ from $X$ onto $Y$ so that 
$$p(S_B(b)x)=S_B(b)p(x)
\tussenen
p(xS_B(b))=p(x)S_B(b)$$
for all $b\in B$.

\snl\bf Proof\rm:
Let $q$ be a linear projection map from $X$ onto $Y$. Define $p:X\to Y$ by
$$p(x)=S_B(E_{(1)})q(E_{(2)}xS_B(E'_{(1)}))E'_{(2)}$$
where $E_{(1)}\ot E_{(2)}$ and $E'_{(1)}\ot E'_{(2)}$ denote two copies of $E$ in the Sweedler notation. The proof is completed using arguments as before. 
\hfill $\square$ \einspr

In the regular case, we get $p(cx)=cp(x)$ and $p(xc)=p(x)c$ for all $c\in C$ while in the semi-regular case, we only get this for $c$ in the subalgebra $S_B(B)$ of $C$.
\snl
Similarly we have the result for $B$-bimodules.

\inspr{4.8} Proposition \rm
Assume that $X$ is a unital $B$-bimodule and $Y$ a submodule. Then there is a projection map $p$ from $X$ onto $Y$ so that 
$$p(S_C(c)x)=S_C(c)p(x)
\tussenen
p(xS_C(c))=p(xS_C(c)$$
for all $c\in C$.

\snl\bf Proof\rm:
Let $q$ be a linear projection map from $X$ onto $Y$. Define $p:X\to Y$ by
$$p(x)=E_{(1)}q(S_C(E_{(2)})xE'_{(1)})S_C(E_{(2)})$$
where $E_{(1)}\ot E_{(2)}$ and $E'_{(1)}\ot E'_{(2)}$ denote two copies of $E$ in the Sweedler notation. The proof is completed as before. 
\hfill $\square$ \einspr

Again, if $E$ is regular, we get here $p(bx)=bp(x)$  and $p(xb)=p(x)b$ for all $b\in B$. If $E$ is semi-regular, we get this only for $b$ in the subalgebra $S_C(C)$ of $B$.
\snl
We see that, at least in the regular case, submodules have complements and so they split into irreducible components. We will use this next to obtain the structure of the base algebras.
\nl
\it Structure of the base algebras \rm
\nl
We will only be able to obtain useful results about the structure of the base algebras $B$ and $C$  when $E$ is regular. Nevertheless, we start this item with the following result.
\inspr{4.9} Proposition \rm
Let $B$ and $C$ be non-degenerate algebras and $E$  any separability idempotent in $M(B\ot C)$. \newline
i) For any element $b\in B$, the space $bB$ is a finite-dimensional right ideal of $B$ and it contains $b$. Similarly, for any element $c\in C$, the space $Cc$ is a finite-dimensional left ideal containing $c$. \newline 
ii) For any element $b\in B$, the space $S_C(C)b$ is a finite-dimensional subspace of $B$ containing $b$. Similarly, for any element $c\in C$, the space $cS_B(B)$ is a finite-dimensional subspace of $C$ containing $c$.

\snl\bf Proof\rm:
i) Take an element $b\in B$. Write $(b\ot 1)E=\sum_i p_i\ot q_i$ with each $p_i$ in $B$ and $q_i\in C$.  
For all $c\in C$ and all linear functionals $\omega$ on $C$ we have $bb'=\sum_i \omega(q_ic)p_i$ where $b'=(\iota\ot\omega(\,\cdot\,c))E$. As all elements in $B$ are of the form $b'$, we find that $bB$ is contained in the subspace of $B$ spanned by the elements $(p_i)$. This proves that $bB$ is a finite-dimensional right ideal in $B$. Because $B$ has local units, we have $b\in bB$. 
\snl
Similarly for $C$.
\snl
ii) Take again $b\in B$. We know that also $((\iota\ot S_C)E)(1\ot b)\in B\ot B$. An argument as above will give that $S_C(C)b$ is finite-dimensional. To show that also here $b$ belongs to this subspace, we use that $S_C$ is non-degenerate so that $b$ is a sum of elements of the form $S_C(c)b'$, together with the existence of local units in $C$. As a consequence, there is an element $e\in C$ so that $S_C(e)b=b$.
\snl
Similarly $cS_B(B)$ is finite-dimensional for all $c\in C$ and $c\in S_B(B)$.
\hfill $\square$ \einspr 

In the semi-regular case, we can already do better.

\iinspr{4.10} Proposition \rm
If $E$ is semi-regular, also $Bb$ is a finite-dimensional left ideal containing $b$ and $cC$ is a finite-dimensional right ideal containing $c$. So, in this case, also $BbB$ will be a finite-dimensional two-sided ideal containing $b$ and $CcC$ a finite-dimensional two-sided ideal containing $c$.

\snl\bf Proof\rm:
i) Take $b\in B$ and now write
$E(b\ot 1)=\sum_i p_i\ot q_i$ with $p_i\in B$ and $q_i\in C$ for all $i$. Again for all $c$ in $C$ and all linear functionals $\omega$ on $C$ we have $b'b=\sum_i \omega(q_ic)p_i$ where $b'=(\iota\ot\omega(\,\cdot\,c))E$. So $Bb$ is finite-dimensional. It contains $b$ as we have local units. Similarly for $C$.
\snl
ii) Take again $b\in B$. From Proposition 4.9.i  we know that $bB$ is finite dimensional. As now also $Bb'$ is finite-dimensional for each $b'\in bB$, we get that $BbB$ is still finite-dimensional. It contains $B$ because we have local units. Similarly for $C$.
\hfill $\square$ \einspr

We will now apply Proposition 4.10 to get the following result.

\iinspr{4.11} Proposition \rm
Assume that $E$ is a regular separability idempotent. Assume that $J$ is a two-sided ideal of $B$. Then it is of the form $Be$ where $e$ is a central idempotent in $M(B)$. If moreover $J$ is finite-dimensional, the element $e$ belongs to $J$ and hence it is the (unique) identity in $J$. A similar result holds for $C$.

\snl\bf Proof\rm:
Let $J$ be a two-sided ideal of $B$. If we consider $B$ as a bi-module over itself, it is unital and $J$ is a submodule. We can apply Proposition 4.8. So there is a projection map $p:B\to J$ such that $p(b'b)=b'p(b)$ and $p(b'b)=p(b')b$. This implies the existence of an element $e\in M(B)$ satisfying $p(b)=eb$ and $p(b)=be$ for all $b$. Hence it is a central element. And because $p$ is a projection map on $J$ we have that $e$ is an idempotent and that $be=b$ when $b\in J$.
\snl
Now we use that $J$ is finite-dimensional. Because $B$ has local units, there exists an element $f\in B$ so that $fb=b$ for all $b\in J$. Then $feb=eb$ for all $b\in B$. This implies $fe=e$ so that $e\in Be$ and $e\in J$. Then $e$ is an identity in $J$. 
\hfill $\square$ \einspr

If the separability idempotent is not assumed to be regular, it is not clear what is still possible. That is even true for the semi-regular case. We refer to Section 5.
\snl
If we combine this result with the property that any element of $B$ belongs to a finite-dimensional two-sided ideal, we arrive at the following structure theorem.

\iinspr{4.12} Proposition \rm
If $E$ is a regular separability idempotent in $M(B\ot C)$, then the algebras $B$ and $C$ are semi-simple. This means that they are direct sums of finite-dimensional matrix algebras over $\Bbb C$.
\snl\bf Proof\rm: 
The proof of this result is standard. 
\snl
i) Take any finite-dimensional two-sided ideal $J$ of $B$ and assume that $e$ is the identity in $J$. It is a central element in $B$. If $J$ is not simple, we have a two-sided ideal $J_1$ of $B$, strictly contained in $J$. It also has an identity $e_1$, central in $B$. Define $J_2=e_2B$ where $e_2=e-e_1$. It follows that $J$ is the direct sum of finitely many simple two-sided ideals. Because we are working with algebras over the field of complex numbers, these components are all full matrix algebras. 
\snl
ii) Because any element of $B$ belongs to a finite-dimensional two-sided ideal, the algebra $B$ will be a direct sum of finite-dimensional simple two-sided ideals.
\hfill $\square$ \einspr

We can also give the following structure theorem of regular separability idempotents.

\iinspr{4.13} Theorem \rm 
Let $B$ and $C$ be non-degenerate algebras and $E$ a regular separability idempotent in $M(B\ot C)$. Then $B$ is the direct sum $\sum_\alpha B_\alpha$ where each $B_\alpha$ is isomorphic with the algebra of $n_\alpha\times n_\alpha$-matrices over $\Bbb C$ for some $n_\alpha\in \Bbb N$. Denote by $e_\alpha$ the identity in $B_\alpha$ and let $f_\alpha=S_B(e_\alpha)$. Then $C$ is the direct sum $\sum_\alpha C_\alpha$ where $C_\alpha=Cf_\alpha$. For each $\alpha$, the element $E_\alpha$ defined in $B_\alpha\ot C_\alpha$ by $E_\alpha=E(e_\alpha\ot f_\alpha)$ is a separability idempotent in $B_\alpha\ot C_\alpha$ and $E$ is the direct product of these components $E_\alpha$ as in Proposition 3.10.

\snl\bf Proof\rm: 
i) We have shown in the previous proposition that $B$ is the direct sum $\sum_\alpha B_\alpha$ where each $B_\alpha$ is isomorphic with the algebra of $n_\alpha\times n_\alpha$-matrices over $\Bbb C$ for some $n_\alpha\in \Bbb N$. The identity $e_\alpha$  in $B_\alpha$ is a central idempotent in $B$. If we let $f_\alpha=S_B(e_\alpha)$, then we get a central idempotent in $C$ because $S_B$ is an anti-isomorphism from $B$ to $C$. Then $C$ is the direct sum $\sum_\alpha C_\alpha$ where $C_\alpha=Cf_\alpha$. Remark that $S_C(f_\alpha)=e_\alpha$ because $e_\alpha$ is central so that $S_CS_B(e_\alpha)=e_\alpha$ (cf.\ a remark after Proposition 2.8 in Section 2). 
\snl
ii) For each $\alpha$ we have $E(e_\alpha\ot 1)=E(1\ot f_\alpha)$. This will give that $E$ only has components in $B_\alpha\ot C_\alpha'$ for $\alpha=\alpha'$. Then it is rather straightforward to see that $E$ is the direct product of these component as in Proposition 3.10 of the previous section.
\hfill $\square$ \einspr

In Section 3, we have seen that any separability idempotent in $M_n(\Bbb C)$ is of the form as in Proposition 3.2 (see Proposition 3.3). This, together with the theorem above, gives a complete description of regular separability idempotents.
\nl
If the algebras $B$ and $C$ are $^*$-algebras and if $E$ is a separability idempotent in $M(B\ot C)$, we know by Proposition 3.7 of Section 3 that it is automatically regular. For any self-adjoint element $b$ in $B$, the ideal $BbB$ is a self-adjoint ideal and hence has the form $Be$ where now $e$ is a self-adjoint central idempotent. It follows that $B$ will be the direct sum as a $^*$-algebra of finite-dimensional matrix algebras with the standard involution.
\snl
It is also a consequence of this property that the distinguished linear functionals $\varphi_B$ and $\varphi_C$ are positive if the algebras are $^*$-algebras and if $E$ is self-adjoint.
\snl
We have also shown in Proposition 3.5 that regularity is automatic in the abelian case. If now $B$ (and hence $C$) is abelian, it means that the components all have to be one-dimensional and it follows that $B$ has the form $K(X)$ for a set $X$ and the separability idempotent $E$ has to be the one of Example 3.6.
\snl
In the general case, and even in the semi-regular case, it seems not obvious how the splitting results for modules, as obtained in the Propositions 4.3 to 4.8, can be used to say something about the structure of the base algebras. We refer to Section 5 for more remarks about this.
\nl\nl

\bf 5. Conclusions and further research \rm
\nl
In this paper, we have studied separability idempotents. 
\snl
We have seen in Section 3 that they arise naturally in the theory of weak multiplier Hopf algebras (as the canonical idempotent), as well as in the case of discrete quantum group (as $\Delta(h)$ for the normalized cointegral, if it exists). The two cases are different in an essential way. And they give rise to two different sets of problems and possible further research. This is what we discuss below.
\snl
The first case suggests the existence of non-regular separability idempotents. Recall that a weak multiplier Hopf algebra is regular if and only if its antipode is bijective. There is a similar result for multiplier Hopf algebras, as well as for ordinary Hopf algebras. Now, there are examples of Hopf algebras with an antipode that is not bijective. So, there are examples of multiplier Hopf algebras and weak multiplier Hopf algebras that are not regular. However, these examples are not good enough to find examples of non-regular separability idempotents. For this we need more exotic cases. We would like e.g.\ to find a weak multiplier Hopf algebra $A$ where the antipode does not map $A$ to $A$ but only to $M(A)$ as suggested in the general theory. Such examples are not yet known. We would also want examples of coproducts that do not satisfy regularity. This means that we would want a coproduct $\Delta$ on an algebra $A$ with the property that elements of the type $(1\ot a)\Delta(b)$ are not necessarily in $A\ot A$. Even such examples are not yet known, but should certainly exist.
\snl
From the general theory of weak multiplier Hopf algebras, one can expect the existence of non-regular separability idempotents, but the lack of non well-behaved examples of non-regular weak multiplier Hopf algebras makes it difficult to find such non-regular separability idempotents.
\snl
Moreover, as we have seen, regularity of a separability idempotent is automatic in various cases: (1) When the algebras are finite-dimensional, (2) when they are abelian and (3) in the involutive case. Also this suggest that it may be hard to find examples of non-regular separability idempotents.
\nl
Let us now turn our attention to the example coming from a discrete quantum group. Take the case of a  multiplier Hopf algebra $(A,\Delta)$ with a cointegral, but so that the left and the right cointegrals are different. In this case, the left cointegral $h$ necessarily has to satisfy $h^2=0$ and so $\Delta(h)$ is not an idempotent. Still, it has a lot of properties in common with separability idempotents. We have e.g.\ that it is full in the sense that the left and right legs are all of the algebra $A$. Moreover, we have $(1\ot a)\Delta(h)=(S(a)\ot 1)\Delta(h)$ where $S$ is the antipode of $A$. For the formula with $S$ on the other side, we would need the right cointegral. However, for the left cointegral we have
$$\Delta(h)(1\ot a)=\sum_{(a)}\Delta(ha_{(1)})(1\ot S(a_{(2)})$$
for all $a$. By the uniqueness of left cointegrals, there exists a homomorphism $\gamma:A\to \Bbb C$ given by $ha=\gamma(a)h$. Therefore we have that $\Delta(h)(1\ot a)=\Delta(h)(1\ot S'(a))$ where $S'(a)=\sum_{(a)}\gamma(a_{(1)})S(a_{(2)})$. We see that even in this case, we still have the antipodal on both sides.
Finally, the left and the right integral $\varphi$ and $\psi$ of the discrete quantum group will give distinguished linear functionals like $\varphi_C$ and $\varphi_B$ in the case of a separability idempotent.
\snl
Recall also that examples in the finite-dimensional case naturally appear where $E^2=0$ instead of $E^2=E$. See a remark before Proposition 3.2 in Section 3.
\snl
Such examples are not possible in the involutive case with self-adjoint elements $E$. We see this in the case of these finite-dimensional examples, as well as in the case of a discrete quantum group. In the situation of a weak multiplier Hopf algebras, where $E$ plays the role of $\Delta(1)$, we will never end up with an example with $E^2=0$.
\snl
Still, it seems worthwhile to develop a more general theory along these lines. 
\nl\nl

\bf References \rm
\nl
{[\bf Dr-VD-Z]} B.\ Drabant, A.\ Van Daele \& Y.\ Zhang: {\it Actions of multiplier Hopf algebra}. Communications in Algebra {\bf 27} (1999), 4117-4172.
\snl
{[\bf Pe]} R.\ Pears: {\it Associative Algebras}. Springer (1982).
\snl
{[\bf VD1]} A.\ Van Daele: {\it Multiplier Hopf algebras}. Trans. Am. Math. Soc.  342(2) (1994), 917-932.
\snl
{[\bf VD2]} A.\ Van Daele: {\it An algebraic framework for group duality}. Adv. in Math. 140 (1998), 323-366.
\snl
{[\bf VD3.v1]} A.\ Van Daele: {\it Separability idempotents and multiplier algebras}. Preprint University of Leuven (2013), arXiv: 1301.4398v1 [math.RA].
\snl
{[\bf VD-W1]} A.\ Van Daele \& S.\ Wang: {\it The Larson-Sweedler theorem for multiplier Hopf algebras}. J.\ of Alg. 296 (2006), 75--95. 
\snl 
{[\bf VD-W2]} A.\ Van Daele \& S.\ Wang: {\it Weak multiplier Hopf algebras. Preliminaries, motivation and basic examples}. Proceedings of  the conference 'Operator Algebras and Quantum Groups (Warsaw, September 2011), Banach Center Publications, volume 98 (2012), 367--415. See also arXiv: 1210.3954 [math.RA] 
\snl 
{[\bf VD-W3]} A.\ Van Daele \& S.\ Wang: {\it Weak multiplier Hopf algebras I. The main theory}. Preprint University of Leuven and Southeast University of Nanjing (2012). To appear in Crelles Journal. Doi: 10.1515/crelle-2013-0053. See also Arxiv: 1210.4395 [math.RA].
\snl
{[\bf VD-W4.v1]} A.\ Van Daele \& S.\ Wang: {\it Weak multiplier Hopf algebras II. The source and target algebras} Preprint University of Leuven and Southeast University of Nanjing (2012), arXiv: 1403.7906v1 [math.RA]
\snl
{[\bf VD-W4.v2]} A.\ Van Daele \& S.\ Wang: {\it Weak multiplier Hopf algebras II. The source and target algebras} Preprint University of Leuven and Southeast University of Nanjing (2012). Arxiv: 1403.7906v2 [math.RA]
\snl
{[\bf VD-Z]} A.\ Van Daele and Y.\ Zhang: {\it Multiplier Hopf algebras of  discrete type}. J. Algebra, 214(1999), 400-417.
\snl
{[\bf Ve]} J.\ Vercruysse: {\it Local units versus local projectivity dualisations: Corings with local structure maps}. Commun. in Alg. 34 (2006) 2079–2103.
\snl
{[\bf W]} Charles A. Weibel: {\it An Introduction to Homological Algebra}. Cambridge University Press (1995). 
\snl

\end